\documentclass[a4paper,12pt]{amsart}
\usepackage{verbatim}

\usepackage[utf8]{inputenc}
\usepackage[T1]{fontenc}
\usepackage{setspace}
\usepackage{indentfirst}

\usepackage{geometry}
\usepackage[hyperfootnotes=false]{hyperref}

\usepackage{times}
\usepackage{amsthm}
\usepackage{amsmath}
\usepackage{amsfonts}
\usepackage{amssymb}
\usepackage{enumitem}
\usepackage{color}
\usepackage{mathrsfs}
\usepackage{stmaryrd}
\SetSymbolFont{stmry}{bold}{U}{stmry}{m}{n}

\usepackage{bm}
\usepackage{mathrsfs}
\usepackage[bbgreekl]{mathbbol}
\DeclareMathSymbol\bbDelta
\mathord{bbold}{"01}  

\usepackage{tikz-cd}
\usepackage{tikz}
\usetikzlibrary{matrix,arrows,decorations.pathmorphing}
\usepackage[all]{xy}
\usepackage{graphicx}
\usepackage{subfigure}

\usepackage[toc,page]{appendix}

\usepackage{todonotes}
\reversemarginpar

\newcommand{\R}{{\mathbb R}}
\newcommand{\C}{{\mathbb C}}
\newcommand{\Z}{{\mathbb Z}}

\newcommand{\Q}{{\mathbb Q}}
\newcommand{\T}{{\mathbb T}}
\newcommand{\G}{{\mathbb G}}
\newcommand{\g}{\mathfrak{g}}

\newcommand{\hhh}{\mathfrak{h}}
\newcommand{\nnn}{\mathfrak{n}}
\newcommand{\n}{^{-1}}
\newcommand{\Hom}{\operatorname{Hom}}

\newcommand{\Id}{\text{Id}}
\newcommand{\Ccal}{{\mathcal C}}
\newcommand{\est}{\text{est}}
\newcommand{\Poly}{\text{Poly}}
\newcommand{\pr}{\operatorname{pr}}
\newcommand{\Fix}{\operatorname{Fix}}

\newcommand{\Mat}{\operatorname{Mat}}

\newtheorem{theorem}{Theorem}[section]
\newtheorem{prop}[theorem]{Proposition}
\newtheorem{cor}[theorem]{Corollary}
\newtheorem{lemma}[theorem]{Lemma}

\theoremstyle{definition}
\newtheorem{definition}[theorem]{Definition}
\theoremstyle{definition}
\newtheorem{remark}[theorem]{Remark}
\theoremstyle{definition}
\newtheorem{example}[theorem]{Example}

\newcommand\xqed[1]{%
  \leavevmode\unskip\penalty9999 \hbox{}\nobreak\hfill
  \quad\hbox{#1}}
\newcommand\tria{\xqed{$\lhd$}}

\begin{document}

\title{Poisson Structures and Potentials}

\author{Anton Alekseev}
\address{Section of Mathematics, University of Geneva, 2-4 rue du Lièvre, c.p. 64, 1211 Genève 4, Switzerland}
\email{Anton.Alekseev@unige.ch}

\author{Arkady Berenstein}
\address{Department of Mathematics, University of Oregon, Eugene, OR 97403, USA}
\email{arkadiy@uoregon.edu}

\author{Benjamin Hoffman}
\address{Department of Mathematics, Cornell University, 310 Malott Hall, Ithaca, NY 14853, USA}
\email{bsh68@cornell.edu}

\author{Yanpeng Li}
\address{Section of Mathematics, University of Geneva, 2-4 rue du Lièvre, c.p. 64, 1211 Genève 4, Switzerland}
\email{yanpeng.li@unige.ch}

\footnotetext{\emph{Keywords:} Poisson structures, Poisson-Lie groups, Potentials, Tropicalization \\
MSC 20G42, 53D17}

\begin{abstract}
   We introduce a notion of weakly log-canonical Poisson structures on positive varieties with potentials. Such a Poisson structure is log-canonical up to terms dominated  by the potential. To a compatible real form of a weakly log-canonical Poisson variety we assign an integrable system on the product of a certain real convex polyhedral cone (the tropicalization of the variety) and a compact torus.
   
  We apply this theory to the dual Poisson-Lie group $G^*$ of a simply-connected semisimple complex Lie group $G$. We define a  positive structure and potential on $G^*$ and show that the natural Poisson-Lie structure on $G^*$ is weakly log-canonical with respect to this positive structure and potential. 
   For $K \subset G$ the compact real form, we show that the real form $K^* \subset G^*$ is  compatible  and prove that the corresponding integrable system is defined on the product of the decorated string cone and the compact torus of dimension $\frac{1}{2}({\rm dim} \, G - {\rm rank} \, G)$.
\end{abstract}

\maketitle

\hskip 9cm
{\em To the memory of Bertram Kostant}

\vskip 0.5cm

\section{Introduction}


Varieties with potentials play a crucial role in mirror symmetry and mathematical physics.
A potential is a rational function on an algebraic variety. We are interested in algebraic varieties with positive structures. That is, an open embedding of a split algebraic torus to the variety so as the potential pulls back to a subtraction-free rational function. Assume for simplicity that the potential can be written as a Laurent polynomial
$$
\Phi(x_1, \dots, x_m) = \sum_{\mathbf{a} \in \mathbb{Z}^n\backslash \{0\}} c_\mathbf{a} x^\mathbf{a},
$$
where $\mathbf{a}=(a_1, \dots, a_m)$ and $c_\mathbf{a} \geq 0$. To such a potential one can assign a rational convex polyhedral cone $\mathcal{C}_\Phi \subset \mathbb{Z}^m$ given by inequalities
$$
\sum_{k=1}^m a_k \xi_k <  0 \,\, {\rm for} \,\, {\rm all} \,\, \mathbf{a} \,\, {\rm with} \,\, c_\mathbf{a} >0,
$$
where the $\xi_i$ are coordinates on $\Z^m$. One can also define this cone using the tropicalized potential $\Phi^t: \mathbb{Z}^m \to \mathbb{Z}$ given by 
$$
\Phi^t(\xi) = \max_{\mathbf{a}; \, c_\mathbf{a} >0} \, \sum_{k=1}^m a_k \xi_k.
$$
Then, $\mathcal{C}_\Phi=\{ \xi \in \mathbb{Z}^m; \Phi^t(\xi) < 0\}$. 
We will denote by $\mathcal{C}_\Phi(\mathbb{R})$ 
the open convex polyhedral cone
$$
\Ccal_{\Phi}(\R)=\{ \xi \in \mathbb{R}^m; \Phi^t(\xi) < 0\}\subset \R^m
$$
defined by the same inequalities. 
 We say that a positive function $f$ on $X$ is dominated by the potential $\Phi$ if $f^t(\xi) < 0$ for all $\xi \in \mathcal{C}_\Phi$. Positive functions dominated by the potential $\Phi$ form a semiring under addition and multiplication. If a rational function $f$ can be decomposed as $f=f_+-f_-$, where $f_+$ and $f_-$ are both dominated by $\Phi$, then we say $f$ is weakly dominated by $\Phi$. Rational functions weakly dominated by $\Phi$ form a ring under addition and multiplication.

Important examples of positive varieties with potentials are complete and partial flag varieties for semisimple algebraic groups \cite{BK}. An example important for us is the Borel subgroup $B \subset G$.
For instance, for $G={\rm SL}_2$ elements of $B$ are upper triangular matrices
\[
b =
\left(
\begin{array}{ll}
b_{11} & b_{12} \\
0 & b_{11}^{-1}
\end{array}
\right) .
\]
In this case, the positive structure is defined by the matrix entries of $b$, and the potential \cite{BK1,Ri} is given by
$$
\Phi(b) = \frac{b_{11} + b_{11}^{-1}}{b_{12}}.
$$
The potential cone is the following cone in $\mathbb{Z}^2$:
\begin{equation} \label{GZ_2}
\mathcal{C}_\Phi = \{ (\xi_{11}, \xi_{12}) \in \mathbb{Z}^2; \xi_{12} > \xi_{11} > - \xi_{12} \}.
\end{equation}
One can show that regular functions dominated by the potential are all non-zero Laurent polynomials of the form
$$
f(b_{11}, b_{12}) = \sum_{k \geq 1} \sum_{l=-k}^k c_{kl} b_{11}^l b_{12}^{-k}
$$
with $c_{kl} \geq 0$.
Note that the cone \eqref{GZ_2} is the $n=2$ instance of the Gelfand-Zeitlin cone \cite{AD, AM} for ${\rm SL}_n$. In fact, for any $n$ there is a positive stucture and the potential on $B \subset {\rm SL}_n$ with potential cone equal to the corresponding Gelfand-Zeitlin cone. 

Varieties arising in applications of mirror symmetry are often symplectic or Poisson. A Poisson structure is a bivector which induces a Poisson bracket on the ring of regular functions on the variety. The Poisson bracket on a positive variety is called log-canonical if 
$$
\{ x_i, x_j \} = c_{ij} x_ix_j,
$$
where $x_1, \dots, x_m$ are toric coordinates (defined by the positive structure) and $c_{ij}$ is a constant matrix. Important examples of Poisson varieties with log-canonical Poisson structures are cluster varieties, and in particular double Bruhat cells including complete and partial flag varieties \cite{KZ}. For instance, in the example of $B \subset {\rm SL}_2$, the log-canonical Poisson bracket is of the form
$$
\{ b_{11}, b_{12} \} = b_{11} b_{12}.
$$
It turns out that this Poisson structure is multiplicative and gives $B$ the structure of a Poisson algebraic group.

The condition of a Poisson structure to be log-canonical is very restrictive. On a positive variety with potential we can generalize it to a notion of \emph{weakly log-canonical} Poisson structures. This means that the Poisson bracket is given by the formulas
$$
\{ x_i, x_j\} = x_ix_j(c_{ij}+ f_{ij}(x)),
$$
where $f_{ij}(x)$ are functions weakly dominated by the potential $\Phi$. 

Given a weakly log-canonical Poisson structure $\pi_X$ on a smooth complex variety $X$, 
consider the real form  $(K,\pi_K) \subset (X,\pi_X)$ defined by the equations $x_i \in \mathbb{R}$ in the toric chart.
Then, to such a structure we assign a constant Poisson bracket on the space
\begin{equation} \label{partial}
\Ccal \times \T,
\end{equation}
where $\Ccal$ is a subcone of $\Ccal_\Phi(\R)$ and $\T\cong (S^1)^r$ is a compact torus of dimension $r$, where $2r$ is the maximal rank of $\pi_K$.
This Poisson bracket has the form
\begin{equation} \label{action_angle}
\{ \xi_i, \xi_j\} = 0, \hskip 0.3cm \{ \phi_i, \phi_j\} =0, \hskip 0.3cm
\{ \xi_i, \phi_j\} = d_{ij},
\end{equation}
where $d_{ij}\in \R$ is determined by the log-canonical part $c_{ij}$ of the weakly log-canonical bracket $\pi_X$.
Here $\xi_i$'s are coordinates on the cone $\Ccal(\mathbb{R})$ and $\phi_i$'s are coordinates on the torus $\T$. We refer to the space \eqref{partial} as a {\em partial tropicalization} of the weakly log-canonical Poisson variety. Up to a change of variables, the Poisson bracket \eqref{action_angle} defines an integrable system on the partial tropicalization.  

Our prime example is the dual Poisson-Lie group $G^*$ of a semisimple complex Lie group $G$ endowed with the standard Poisson-Lie structure. For $G={\rm SL}_2(\C)$, the group $G^*$ is of the form
$$
G^*=\{ (b^+, b^-) \in B \times B_-; b^+_{11}b^-_{11}=1\}.
$$
The positive structure is defined by the following parametrization:
$$
b^+ =
\left(
\begin{array}{ll}
b_{11} & b_{12} \\
0 & b_{11}^{-1}
\end{array}
\right) ,
\hskip 0.5cm
(b^-)^{-1} =
\left(
\begin{array}{ll}
b_{11} & 0 \\
b_{21} & b_{11}^{-1}
\end{array}
\right) ,
$$
and the potential is given by
$$
\Phi_{G^*}= \frac{b_{11} + b_{11}^{-1}}{b_{12}} + \frac{b_{11} + b_{11}^{-1}}{b_{21}}.
$$

To $K \subset G$ the compact real form of $G$ one associates the real form $K^* \subset G^*$. In the case of $K ={\rm SU}(2) \subset {\rm SL}_2(\mathbb{C})$, this real form is defined by equations $b_{11} \in \mathbb{R}_{>0}$ and $b_{21}=\overline{b_{12}}$.
The canonical real Poisson bracket on $K^*$ is given by \cite{LW,STS}:
$$
\{ b_{11}, b_{12}\} = i b_{11}b_{12}, \hskip 0.3cm \{b_{11}, b_{21}\} = -i b_{11} b_{21}, \hskip 0.3cm \{ b_{12}, b_{21}\} = i (b_{11}^2 - b_{11}^{-2}).
$$
Note that the first two expressions are log-canonical on the nose whereas the third expression has no log-canonical part and the corresponding function $f(x)$ is of the form
$$
i \left( \frac{b_{11}^2}{b_{12} b_{21}} -  \frac{b_{11}^{-2}}{b_{12} b_{21}}\right).
$$
This expression is weakly dominated by the potential $\Phi_{G^*}$. The corresponding partial tropicalization is the product of the Gelfand-Zeitlin cone and the circle $S^1$:
$$
\{ (\xi_{11}, \xi_{12} ) \in \mathbb{R}^2; \xi_{12} > \xi_{11} > - \xi_{12}\} \times S^1
$$
with Poisson bracket
$$
\{ \xi_{11}, \xi_{12}\} =0, \hskip 0.3cm \{ \xi_{11}, \phi\} =1, \hskip 0.3cm
\{\xi_{12}, \phi\} =0,
$$
which is the $n=2$ Gelfand-Zeitlin integrable system.

We will now describe our main results.  Let $G$ be a semisimple complex Lie group with the standard Poisson-Lie structure and $G^*$  the dual Poisson-Lie group. The standard positive structures on $B$ and $B_-$ and the respective potentials $\Phi_+$ and $\Phi_-$ give rise to a positive structure $\theta_{G^*}$ and a potential $\Phi_{G^*}$ on $G^*$. 
We prove the following theorem:
\begin{theorem} \label{I1}
In the coordinates given by $\theta_{G^*}$, the standard Poisson-Lie structure on $G^*$ is weakly log-canonical with respect to the potential $\Phi_{G^*}$.
\end{theorem}

Furthermore, for the compact real form $K \subset G$ we consider the partial tropicalization of the real form $K^* \subset G^*$. Our next result is as follows:
\begin{theorem} \label{I2}
The partial tropicalization of $K^*$ is of the form $\mathcal{C} \times \T$, where $\mathcal{C}$ is the extended string cone and $\T$ is a real torus of dimension $d={\rm dim}(N)$. This space is equipped with an integrable system.
\end{theorem}

In the case of $G={\rm SL}_n(\mathbb{C})$, these results (in one of the toric charts) were obtained in \cite{AD}. We establish them in full generality: for all equivalent positive structures on ${\rm SL}^*_n(\mathbb{C})$ and for all semisimple complex Lie groups.

In the next several paragraphs we describe the future directions which give an additional context and motivation for our results:

Dual Poisson-Lie groups $K^*$ of compact connected Poisson-Lie groups are Poisson manifolds with very special properties. In particular, the Ginzburg-Weinstein isomorphism theorem \cite{GW} states that $(K^*, \pi_{K^*})$ is isomorphic as a Poisson manifold to the dual of the Lie algebra $\mathfrak{k}^*$ with the linear Kirillov-Kostant-Souriau (KKS) Poisson bracket $\pi_{\rm KKS}$. Since the bracket $\pi_{\rm KKS}$ is linear, the scaling map $D_s: x \mapsto s x$ is a Poisson map $D_s: (\mathfrak{k}^*, \pi_{\rm KKS}) \to (\mathfrak{k}^*,  s \pi_{\rm KKS})$. Following \cite{Anton, AM}, this observation gives rise to a family of Ginzburg-Weinstein maps 
$$
{\rm gw}_s: (\mathfrak{k}^*, \pi_{\rm KKS}) \to (K^*, s\pi_{K^*}).
$$
We would like to think of the partial tropicalization $\mathcal{C} \times T$ of $K^*$ as of the $s \to \infty$ limit of the Poisson space $(K^*, s\pi_{K^*})$. 

In the forthcoming paper we plan to give the following evidence in support of this idea: the symplectic leaves in $\mathcal{C} \times T$ are labeled by the elements $\lambda \in \overset{\circ}{W}_+$ of the interior of the positive Weyl chamber and their symplectic volumes coincide with symplectic volumes of the corresponding coadjoint orbits $\mathcal{O}_\lambda \subset \mathfrak{k}^*$. The proof makes use of the relation between the Kashiwara crystals for finite dimensional modules of $G$ and of the Langlands dual group $G^\vee$.

A more ambitious project is to define the $s \to \infty$ limit of the Ginzburg-Weinstein map ${\rm gw}_s$. This will be a Poisson map 
$$
{\rm gw}_\infty: (\mathfrak{k}^*, \pi_{\rm KKS}) \to (\mathcal{C} \times T, \pi_\infty),
$$
where $\pi_\infty$ is the constant Poisson bracket on the partical tropicalization. Among other things, such a map will define interesting completely integrable systems on the top-dimensional coadjoint orbits $\mathcal{O}_\lambda$. For the case of $G={\rm SL}_n(\mathbb{C})$, the map ${\rm gw}_\infty$ can be constructed by combining results of \cite{AM} and \cite{AD}. For the general case of $G$ semisimple this question remains open.

The structure of the paper is as follows: in Section 2 we review the theory of positive varieties, in Section 3 introduce potentials and tropicalizations, in Section 4 we apply this theory to double Bruhat cells and prove the main technical results on functions dominated by the potentials, in Section 5 we discuss weakly log-canonical Poisson brackets and prove Theorem \ref{I1}. In section 6 we consider real forms of positive varieties, define partial tropicalization, and prove Theorem \ref{I2}.

{\bf Acknowledgements.} We are grateful to I. Davydenkova who shared with us her unpublished results on the case of $G={\rm Sp}(4)$, to M. Podkopaeva who participated in the initial stage of this project, to G. Koshevoy for his useful comments on the string cone, and to J. Lu for her helpful comments on an earlier draft. This work became possible thanks to the NCCR SwissMAP, which organized a series of meetings in Les Diablerets where A.B. participated and the Master Class in Geneva where B.H. took part during the academic year 2016-2017. A.B. expresses his gratitude for hospitality and support during his visits to Geneva in 2016 and in 2017. A.A. and Y.L. were supported in part by the ERC project MODFLAT and by the grants 165666 and 159581 of the Swiss National Science Foundation. B.H. was supported by the National Science Foundation Graduate Research Fellowship under Grant Number DGE-1650441.

\section{Positivity Theory}

In this Section we review the theory of positive varieties and positive maps. One of the central notions is tropicalization, which assigns a lattice to a positive variety. We mainly follow the work of Berenstein and Kazhdan, see \cite{BK1} and \cite{BK}.

\subsection{Algebraic Tori and Positive Maps}

Let $\G_a$ be the additive group and $\G_m$ the multiplicative group defined over $\mathbb{Q}$. By definition, the coordinate algebra $\mathbb{Q}[\G_a]$ is the polynomials over $\mathbb{Q}$ in one variable $\mathbb{Q}[x]$ and the coordinate algebra $\mathbb{Q}[\G_m]$ is the Laurent polynomials $\mathbb{Q}[x, x^{-1}]$.

Consider a split algebraic torus $S\cong \G_m^n$. Denote by $S_t={\rm Hom}(S, \G_m)$ the character lattice of $S$ and by $S^t = {\rm Hom}(\G_m, S)$ the cocharacter lattice. The lattices $S_t$ and $S^t$ are in natural duality, in particular $S_t \cong {\rm Hom}(S^t, \mathbb{Z})$. The coordinate ring $\mathbb{Q}[S]$ is the group algebra (over $\mathbb{Q}$) of the lattice $S_t$. That is, every element $f \in \mathbb{Q}[S]$ can be written as
\begin{equation} \label{chicchi}
f= \sum_{\chi \in S_t} c_\chi \chi,
\end{equation}
where only a finite number of coefficients $c_\chi$ are non-zero.

 Denote by $P[S]$ the set of non-zero elements of the form \eqref{chicchi} in which all $c_\chi$ are non-negative rational numbers.
 By construction, $P[S] \subset \mathbb{Q}[S]$ is a sub-semiring without zero.

Denote by $\mathbb{Q}(S)$ the field of rational functions on $S$, and by $P(S) \subset \mathbb{Q}(S)$ the subset of elements of $\mathbb{Q}(S)$ represented as fractions $f/g$ with $f,g \in P[S]$. By construction, $P(S)$ is a semifield without zero.

\begin{example} \label{trivial}
Note that the polynomial $x^2 -x +1=(x^3+1)/(x+1)$ belongs to $P(\G_m)$ but does not belong to $P[\G_m]$, and the polynomial $x^2-2x+1=(x-1)^2$ belongs neither to $P[\G_m]$ nor to $P(\G_m)$.  
\end{example}
%

\begin{definition}
A rational map $\phi:S\to S'$ between split algebraic tori is positive if
for every character $\chi: S' \to \G_m$ the composition $\chi \circ \phi \in P(S)$ is a positive rational function on $S$.
\tria \end{definition}

\begin{example}\label{coordchange}
  Let $S=S'=\G_m^3$. The rational map $\phi:S\to S'$ defined by
  \[
    \phi(x_1,x_2,x_3)=\left(\frac{x_2x_3}{x_1+x_3},x_1+x_3,\frac{x_1x_2}{x_1+x_3}
    \right)
  \]
  is positive.
\end{example}

\begin{prop}
  A rational map $\phi:S\to S'$ is positive if and only if the pullback $\phi^*$ restricts to a semifield homomorphism on positive rational functions:
  \[
    \phi^*:P(S')\to P(S).
  \]
\end{prop}
\begin{proof}
  Let $y_1,\dots,y_m\in P(S')$ be the standard basis of the character lattice $S'_t$ of $S'$ given by the splitting of $S'$, and let
  $\phi_i = y_i \circ \phi \in \mathbb{Q}(S)$ the components of $\phi$. Assume $\phi$ is positive, then $\phi_i\in P(S)$ for each $i$. We show that $\phi^* f\in P(S)$ whenever $f\in P(S')$. Indeed, $\phi^* f(y_1,\dots,y_m)=f(\phi_1,\dots,\phi_m)\in P(S)$ is a subtraction-free rational polynomial in $\phi_1,\dots,\phi_m$.

  Conversely, assume $\phi^* f\in P(S)$ whenever $f\in P(S')$. Then for each basic character $y_i$ on $S'$, we have $\phi^*y_i=y_i \circ \phi = \phi_i \in P(S)$.
\end{proof}

Denote by $\mathbf{PosTori}$ the category with objects split algebraic tori and arrows positive rational maps. The previous proposition shows $P(\cdot)$ defines a functor from $\mathbf{PosTori}^{op}$ to the category of semifields. However, the next example shows the situation is not so straight forward when we take positive regular functions $P[S]$ of a split algebraic torus $S$; see also Remark \ref{posregnonoremark} below.

\begin{example} \label{posregnono}
Define a rational map $F:\G_m^2\to \G_m^2$ by $F(t_1,t_2)=\left(\frac{t_1t_2}{t_1+t_2},\frac{t_2^2}{t_1+t_2}\right)$. Clearly, $F$ is invertible and $F^{-1}$ is given by $F^{-1}(t_1,t_2)=\left(\frac{t_1}{t_2}\cdot(t_1+t_2),t_1+t_2\right)$. In particular, $F$ is an isomorphism in the category $\mathbf{PosTori}$. Let $f:=\frac{(t_1^3+t_2^3)(t_1+t_2)}{t_2^2}\in P[\G_m^2]$. Then 
\[
f\circ F=\frac{\left(\left(\frac{t_1t_2}{t_1+t_2}\right)^3+\left(\frac{t_2^2}{t_1+t_2}\right)^3\right)\left(\frac{t_1t_2}{t_1+t_2}+\frac{t_2^2}{t_1+t_2}\right)}{\left(\frac{t_2^2}{t_1+t_2}\right)^2}
=\frac{t_1^3+t_2^3}{t_1+t_2}=t_1^2-t_1t_2+t_2^2\ .
\]
Thus, $f\circ F\notin P[\G_m^2]$.
\end{example}

\subsection{Tropicalization of positive maps}

Following \cite{BK}, 
to each positive rational map $\phi: S\to S'$ we will associate a \emph{tropicalized map} $\phi^t:S^t\to (S')^t$ as follows.

{\em Step 1.} Let $\phi: S \to \G_m$ be a positive regular function, that is $\phi \in P[S]$ and it admits the form \eqref{chicchi} with $c_\chi \geq 0$. Define $\phi^t: S^t \to \G_m^t = \mathbb{Z}$ by formula
$$
\phi^t(\xi) = \max_{\chi; \, c_\chi >0} \langle \chi, \xi \rangle,
$$
where $\langle \cdot, \cdot \rangle: S_t  \times S^t \to \mathbb{Z}$ is the canonical pairing.

\vskip 0.2cm

{\em Step 2.} Let $\phi: S \to \G_m$ be a positive rational function. That is, $\phi=f/g$ with $f,g \in P[S]$, then 
$$
\phi^t=f^t - g^t.
$$
One can show that the right hand side is independent of the presentation of the positive fraction. Note that the assignment $ \phi \mapsto \phi^t$ is a homomorphism of semifields mapping $P(S)$ to the semifield of $\mathbb{Z}$-valued functions on $S^t$. This semifield has operations of tropical addition $(f, g) \mapsto \max\{f,g\}$ and tropical multiplication $(f,g) \mapsto f+g$.

\begin{example}
 Consider $S=\G_m$ and $\phi=(x^3+1)/(x+1)$. Then, $\phi^t(\xi)=\max(3\xi, 0) - \max(\xi, 0)= 2 \max(\xi, 0)$. Note that for any $a,b,c,d \in \mathbb{Q}_{>0}$ the function $\phi'=(ax^3+b)/(cx+d)$ has the same tropicalization $(\phi')^t=\phi^t$.
\end{example}

{\em Step 3.} Let $\phi: S \to S'$ be a positive rational map. Define $\phi^t: S^t \to (S')^t$ as the unique map such that for every character $\chi \in S'_t$ and for every cocharacter $\xi \in S^t$ we have
$$
\langle \chi, \phi^t(\xi) \rangle = (\chi \circ \phi)^t(\xi).
$$
Let $\phi_1,\dots,\phi_m$ be the components of $\phi$ given by the splitting $S'\cong (\mathbb{G}_m)^m$. Then, in the induced coordinates on $(S')^t,$ we have \[ \phi^t=(\phi_1^t,\dots,\phi_m^t).\]

\begin{example}\label{tropeg}
  The positive rational map from Example \ref{coordchange}
  \begin{align*}
    \phi:\G_m^3 \to \G_m^3 \ :\ 
    (x_1,x_2,x_3)& \mapsto \left(\frac{x_2x_3}{x_1+x_3},x_1+x_3,\frac{x_1x_2}{x_1+x_3} \right)
  \end{align*}
  has tropicalization
  \begin{align*}
    \phi^t:(\G_m^3)^t\cong \Z^3 & \to (\G_m^3)^t\cong \Z^3; \\
    (\xi_1,\xi_2,\xi_3) &\mapsto (\xi_2+\xi_3-\max\{\xi_1,\xi_3\},\max\{\xi_1,\xi_3\}, \xi_1+\xi_2-\max\{\xi_1,\xi_3\}).
  \end{align*}
  Note that $\phi^t$ is linear on the chambers $\xi_1<\xi_3$ and $\xi_1>\xi_3$.
\end{example}
 
Recall that $\mathbf{PosTori}$ the category with objects split algebraic tori and arrows positive rational maps. Denote by $\mathbf{PLSpaces}$ the category with objects finite rank lattices and arrows piecewise $\mathbb{Z}$-linear maps. The following is shown in Section 2.4 of \cite{BK1}:

\begin{prop}\label{functorial}
  The assignment $S \mapsto S^t, \phi \mapsto \phi^t$ defines a functor from the category $\mathbf{PosTori}$ to the category $\mathbf{PLSpaces}$.
\end{prop}

Note that all piecewise linear maps $\phi^t$ obtained by applying the tropicalization functor are homogeneous in the following sense:
\begin{equation}
\label{homogeq}
\phi^t(n \xi) = n \phi^t(\xi)
\end{equation}
for every $n \in \mathbb{Z}_{\geq 0}$. In particular, $\phi^t(0)=0$. 
 
\subsection{Positive Varieties}
 
\begin{definition}
  Let $X$ be an irreducible variety over $\mathbb{Q}$. A \emph{toric chart} is an open embedding  $\theta: S \to X$ from a split algebraic torus $S$ to $X$.
\tria \end{definition}

Since $\theta$ is an open map, it induces an inclusion of coordinate rings. Hence, the coordinate ring of $X$ identifies with a subalgebra of $\mathbb{Q}[S]$. 
   
\begin{example}
  For $X=S$, the identity map $\Id: S\to S$ is a toric chart.
\end{example}

\begin{example}
  The inclusion $\G_m\to \G_a$ is a toric chart on $\G_a$. The corresponding homomorphism of coordinate algebras is the natural embedding $\mathbb{Q}[x] \to \mathbb{Q}[x, x^{-1}]$.
\end{example}

\begin{example}
  Let $X$ be a toric variety, and let $\theta: S\hookrightarrow X$ 
  pick out the open torus of $X$. Then $\theta$ is a toric chart.
\end{example}

\begin{example} \label{charteg1}
  Let $N\subset {\rm SL}_3$ be the group of unipotent upper-triangular matrices. Define $\theta:S=(\G_m)^3\to N$ by
  \begin{align*}
    \theta(x_1,x_2,x_3) & = 
      \begin{pmatrix}
        1 & x_1 & 0 \\
        0 & 1 & 0 \\
        0 & 0 & 1
      \end{pmatrix}
      \begin{pmatrix}
        1 & 0 & 0 \\
        0 & 1 & x_2  \\
        0 & 0 & 1
      \end{pmatrix}
      \begin{pmatrix}
        1 & x_3 & 0 \\
        0 & 1 & 0 \\
        0 & 0 & 1
      \end{pmatrix}=
      \begin{pmatrix}
        1 & x_1+x_3 & x_1x_2 \\
        0 & 1 & x_2 \\
        0 & 0 & 1
      \end{pmatrix}.
  \end{align*}
  Then $\theta$ is a toric chart on $N$.
\end{example}

\begin{example}\label{charteg2}
  Let $N\subset {\rm SL}_3$ be as in the previous example.
  Define $\theta':(\G_m)^3\to N$ by
  \begin{align*}
    \theta'(y_1,y_2,y_3) & = 
      \begin{pmatrix}
        1 & 0 & 0 \\
        0 & 1 & y_1 \\
        0 & 0 & 1
      \end{pmatrix}
      \begin{pmatrix}
        1 & y_2 & 0 \\
        0 & 1 & 0  \\
        0 & 0 & 1
      \end{pmatrix}
      \begin{pmatrix}
        1 & 0 & 0 \\
        0 & 1 & y_3 \\
        0 & 0 & 1
      \end{pmatrix}=
      \begin{pmatrix}
        1 & y_2 & y_2y_3 \\
        0 & 1 & y_1+y_3 \\
        0 & 0 & 1
      \end{pmatrix}.
  \end{align*}
  Then $\theta'$ is a toric chart on $N$.
\end{example}
 
\begin{example}
  Let $A$ be an upper cluster algebra, let $X=\operatorname{Spec} A$, and let $x_1,\dots,x_k$ be a cluster. By the Laurent phenomenon, $A\subset \mathbb{Q}[x_1^\pm,\dots,x_k^\pm]$, and the corresponding map of varieties
  \[
    \G_m^{k} = \operatorname{Spec} \mathbb{Q}[x_1^\pm,\dots,x_k^\pm] \to X
  \]
   is a toric chart. More details can be found in \cite{BFZ}.
\end{example}
 
\begin{definition}
  Let $\theta: S\to X$ and $\theta': S'\to X$ be toric charts on an irreducible variety $X$. If $(\theta)\n\circ\theta': S' \to S$ and $(\theta')\n\circ\theta: S\to S'$ are positive rational maps, then $\theta$ and $\theta'$ are \emph{positively equivalent} toric charts. We define a \emph{positive variety} to be a pair $(X,\Theta_X)$, where $\Theta_X$ is a positive equivalence class of toric charts. If $\theta \in \Theta_X$, we sometimes write $\Theta_X=[\theta]$.
\tria \end{definition}
 
\begin{example}\label{transitionmapex}
  Let $N\subset {\rm SL}_3$ be upper triangular unipotent matrices.
  Consider the toric charts $\theta,\theta'$ on $N$ from Examples \ref{charteg1} and \ref{charteg2}. Then $(\theta')\n\circ \theta$ is the positive rational map $\phi:\G_m^3\to \G_m^3$ from Example \ref{coordchange}. Its inverse has a similar positive expression, so $\theta$ and $\theta'$ are positively equivalent and $[\theta]=[\theta']$.
\end{example}
 
\begin{remark}
  If a birational map $\phi: S\to S'$ is positive, its inverse is not necessarily positive. For instance, the birational map $\G_m \to \G_m$ given by $x\mapsto x+1$ does not have a positive inverse.
\end{remark}

\begin{definition}
  A \emph{positive map} of positive varieties $\phi:(X,\Theta_X)\to (Y,\Theta_Y)$ is a rational map $\phi:X\to Y$ so that for some (equivalently any) $\theta_X\in \Theta_X$ and $\theta_Y\in \Theta_Y$, the rational map $\theta_Y\n \circ\phi\circ \theta_X: S\to S'$ is positive.  We denote by $\mathbf{PosVar}$ the category of positive varieties over $\mathbb{Q}$, with positive rational maps as arrows.
\tria \end{definition}
 
In particular, setting $Y=\G_m$, we see that a rational function $f\in \mathbb{Q}(X)$ is positive if $f\circ\theta_X: S \to \G_m$ belongs to $P(S)$. Denote the set of positive rational functions on $(X,\Theta_X)$ by $P(X,\Theta_X)$. It is clear that $P(X, \Theta_X)$ is a semifield.

\begin{remark} \label{posregnonoremark}
Example \ref{posregnono} shows that we cannot similarly define a ``positive regular semiring,'' $P[X,\Theta_X]$ as the definition will not be independent of the choice of toric chart. However, we indicate in Remark \ref{totally_positive} that in some situations there is an alternate approach.
\end{remark}

 
\begin{definition}
  Let $\mathbf{PosVar}^\bullet$ be the category with objects \emph{framed positive varieties} $(X,\theta)$, where $\theta$ is a toric chart on $X$. An arrow from $(X, \theta_X)$ to $(Y, \theta_Y)$ is a rational map $\phi: X \to Y$ such that 
  $\theta_Y^{-1} \circ \phi \circ \theta_X: S_X \to S_Y$ is a positive map of tori. 
\tria \end{definition}

Let
\[
  F:\mathbf{PosVar}^\bullet\to \mathbf{PosVar}
\]
be the forgetful functor sending $(X,\theta)\mapsto (X,[\theta])$, then it is shown in Claim 3.17 in \cite{BK} that $F$ is an equivalence of categories. We may define an adjoint equivalence 
\[
  G: \mathbf{PosVar}\to \mathbf{PosVar}^\bullet
\]
by simultaneously choosing a representative $\theta\in \Theta_X$ for each positive variety $(X,\Theta_X)$. In fact, any adjoint equivalence to $F$ arises this way.
 
Tropicalization extends in the obvious way to framed positive varieties: if $(X,\theta)$ is a framed positive variety, with $\theta: S\to X$ a toric chart, set
\[
  (X,\theta)^t:=\Hom(\G_m,S)=S^t.
\]
A positive rational map $\phi:(X,\theta_X)\to(Y,\theta_Y)$ has a tropicalization $\phi^t:=(\theta_Y\n\circ \phi\circ \theta_X)^t:(X,\theta)^t\to (Y,\theta)^t$, and by Proposition \ref{functorial}, tropicalization respects composition of positive rational maps. We then have:
 
\begin{prop}\label{tropfunctor}
  Tropicalization defines a functor from the category $\mathbf{PosVar}^\bullet$ to the category $\mathbf{PLSpaces}$.
\end{prop}
 
\begin{example}
  Let $(X, \Theta)$ be a positive variety. Then for every $\theta, \theta' \in \Theta$ the identity map of $X$ induces a map $(X, \theta) \to (X, \theta')$ in $\mathbf{PosVar}^\bullet$. In particular, recall
  Example \ref{transitionmapex}. There, we considered $X=N$ with two toric charts.
  The tropicalization of the induced map was given in 
  Example \ref{tropeg}.
\end{example}
 
In general, if $\theta,\theta':\G_m^k\to X$ are positively equivalent charts, the transition map $(\theta\n\circ \theta')^t:\Z^k\to \Z^k$ is a piecewise $\Z$-linear bijection.
  
Precomposing the functor $(\cdot)^t:\mathbf{PosVar}^\bullet \to \mathbf{PLSpaces}$ with the equivalence $G: \mathbf{PosVar}\to \mathbf{PosVar}^\bullet$ gives a tropicalization of positive varieties. Choosing a different adjoint equivalence $G'$, we see $(G(X,\Theta))^t$ differs from $(G'(X,\Theta))^t$ by a piecewise linear bijection.
 
\section{Potentials}

In this Section we introduce the notion of potentials on positive varieties. Potentials are positive functions, and from our perspective their main role is to define interesting cones on tropicalized varieties.

\begin{definition}
  Let $(X,\theta)$ be a framed positive variety. We distinguish a set of positive rational functions $\Phi:=\{\varphi_1,\dots,\varphi_m\}\subset P(X,[\theta])$, called a \emph{set of potentials on} $(X,\theta)$. The triple $(X,\theta,\Phi)$ is called a \emph{framed positive variety with potential}.

  Since the semiring $P(X,[\theta])$ depends only on the positive structure $[\theta]$ and not the toric chart $\theta$, we define similarly a \emph{positive variety with potentials} as a triple $(X,[\theta],\Phi)$.
  
  We write $\Phi^t=\{\phi_1^t,\dots,\phi_m^t\}$. For $\xi\in (X,\theta)^t$, we write $\Phi^t(\xi)<0$ whenever $\phi_i^t(\xi)<0$ for all $\phi_i^t\in \Phi^t$.
\tria \end{definition}

\begin{definition}
  Let $(X,\theta)$ be a framed positive variety. For $\Phi\subset P(X,[\theta])$, we define the cone
  \begin{align*}
    \Ccal_{\Phi}(X,\theta) & :=\{\xi\in (X,\theta)^t|\Phi^t(\xi)<0\}\subset (X,\theta)^t.
  \end{align*}
  We write $\Ccal_\Phi$ when the variety and positive structure $(X,\theta)$ are evident, and we use $\mathcal{C}_f$ if $\Phi=\{f\}$. If $\Phi$ is distinguished as a potential on $(X,\theta)$, we call $\Ccal_\Phi$ the \emph{(strict) potential cone}.
\tria \end{definition}

\begin{remark}
\label{coneremark}
\begin{enumerate}
    \item 
 In \cite{BK1,BK,BZ}, the authors consider the cone
 \[
    \Ccal_{\Phi}^\le(X,\theta) :=\{\xi\in (X,\theta)^t|\varphi_i^t(\xi)\le 0 \mbox{~for all~} \varphi_i\in \Phi\}\subset (X,\theta)^t
 \] given by non-strict inequalities, whereas we primarily consider the strict potential cone $\Ccal_\Phi$. In Section \ref{lastsection}, we consider the open real cone
 \[
 \Ccal_\Phi(\R)=\{\xi\in (X,\theta)^t\otimes \R | \Phi^t(\xi)<0\}
 \]
 defined by strict inequalities.
 
\item When $\Phi$ restricts to a set of regular functions on the toric chart $\theta:S\hookrightarrow X$, the cone $\Ccal_\Phi(\R)$ is polyhedral and hence convex.

 \item We sometimes consider the function $0$ as a potential on a positive variety $(X,\theta)$. By convention, $0^t=-\infty$ and $\Ccal_0(X,\theta)=(X,\theta)^t$.
 \end{enumerate}
\end{remark}
 
\begin{definition}
  Let $(X,\theta)$ be a framed positive variety and let $\Phi,\Psi\subset P(X,[\theta])$.
  Then $\Psi$ is \emph{dominated} by $\Phi$ if, for all $\xi\in (X,\theta)^t$, if $\Phi^t(\xi)<0$ then $\Psi^t(\xi)<0$. 
\tria \end{definition}

In other words, $\Psi$ is dominated by $\Phi$ if and only if $\Ccal_\Phi(X,\theta) \subset \Ccal_\Psi(X,\theta)$.


\begin{lemma}\label{conesinclusion}
  Let $f:(X,\theta_X)\to  (Y,\theta_Y)$ be a positive rational map. If $\Phi,\Psi\subset P(Y,[\theta_Y])$ and $\Psi$ is dominated by $\Phi$, then, $f^*\Psi$ is dominated by $f^*\Phi$. Here, we write
  \[
    f^*\Phi:=\{f^*\phi_1,\dots,f^*\phi_n\},\text{ where }\Phi=\{\phi_1,\dots,\phi_n\}.
  \]
\end{lemma}

\begin{proof}
  Without loss of generality we assume $\Phi=\{\phi\}$ and $\Psi=\{\psi\}$; see Proposition \ref{singlepotential} below.
  Since tropicalization is a functor, we know $(f^*(\phi))^t=\phi^t \circ f^t$. By assumption, 
  \[
  \phi^t(\xi)<0\Rightarrow \psi^t(\xi)<0\text{ for all }\xi\in (Y,\theta_Y)^t.
  \]
  So,
  \[
  (f^*(\phi))^t=\phi^t( f^t(\xi))<0\Rightarrow \psi^t(f^t(\xi))=(f^*(\psi))^t(\xi)<0\text{ for all }\xi\in (X,\theta_X)^t,
  \]
   and thus $f^*\Psi$ is dominated by $f^*\Phi$.      
\end{proof}

As a consequence of Lemma \ref{conesinclusion}, domination is preserved by positive rational equivalences. 

\begin{definition}
  Let $(X,\Theta_X,\Phi_X)$ and $(Y,\Theta_Y,\Phi_Y)$ be positive varieties with potential. A positive rational map $f:(X,\Theta_X)\to (Y,\Theta_Y)$ is a \emph{map of positive varieties with potential} if for some (equivalently, any) $\theta_X\in \Theta_X$ and $\theta_Y\in \Theta_Y$, we have $f^*\Phi_Y$ is dominated by $\Phi_X$.

  Positive varieties with potential and their maps form a category $\mathbf{PosVarPot}$.
\tria \end{definition}

\begin{definition}
  For $(X,[\theta_X],\Phi)$ a positive variety with potential, let $P_\Phi(X,[\theta_X])$ be the set of positive rational functions on $X$ which are dominated by $\Phi$, which we call the \emph{dominated semi-ring}. It is a semi-ring, with the usual operations $+$ and $\cdot$. Note that $\Ccal_1(X,\theta_X)=\emptyset$, so $P_\Phi(X,[\theta_X])$ has no multiplicative unit in general.

  Let $\tilde{P}_\Phi (X,[\theta_X])\subset \mathbb{Q}(X)$ be the ring (with the usual operations) generated by $P_\Phi(X,[\theta_X])$, which we call the \emph{ring of weakly dominated functions}.
\tria \end{definition}

\begin{prop}
  The assignments $P_\bullet$ and $\tilde{P}_\bullet$ define functors from $\mathbf{PosVarPot}^{op}$ to the category of semirings $\mathbf{SemiRing}$ and $\Q$-algebras $\mathbf{Alg}_\Q$, respectively.
\end{prop}

\begin{proof}
  Follows from Lemma \ref{conesinclusion}.
\end{proof}

\begin{example}\label{potentialeg}
  Consider the framed positive variety $(N,\theta)$ from Example \ref{charteg1}. Define the set of potentials: \[
    \Phi:=\left\{\frac{1}{x_1},\frac{1}{x_3},\frac{x_1+x_3}{x_1x_2},\frac{x_1+x_3}{x_2x_3}\right\}.
  \]
  Tropically, we have
  \[
    \Phi^{t}=\left\{ -\xi_1,-\xi_3, \max\{\xi_1,\xi_3\}-(\xi_1+\xi_2), \max\{\xi_1,\xi_3\}-(\xi_2+\xi_3) \right\}.
  \]
  The potential cone is
  \[
    \Ccal_\Phi=\left\{\xi_1,\xi_2,\xi_3\in \Z^3| \xi_1>0,~\xi_3>0,~\max\{\xi_1,\xi_3\}<\xi_1+\xi_2,~\max\{\xi_1,\xi_3\}<\xi_2+\xi_3 \right\}.
  \]
  Let $f:=\frac{1}{x_1x_2}$, then  \[\Ccal_f=\{\xi_1,\xi_2,\xi_3\in \Z^3| \xi_1+\xi_2>0 \}.\] 
  The function $f$ is dominated by $\Phi$: if a point is in $\Ccal_\Phi$, it is in $\Ccal_f$.
\end{example}

\begin{definition}
  Let $(X,\Theta,\Phi)$ be a positive variety with potential. Let
  \[
    P^{\est}_\Phi(X,\Theta):=\{f|\exists g\in P(X,\Theta)\mbox{~such that~} f+g\in \Poly^+\Phi\},
  \]
  where $\Poly^+\Phi$ is the semi-ring of  polynomials in $\Phi=\{\phi_1,\dots,\phi_n\}$ with positive rational coefficients, and no constant terms. We call  $P^{\est}_\Phi(X,\Theta)$ the semi-ring of \emph{estimate-dominated} functions. We define the \emph{ring of weakly estimate-dominated functions} $\tilde{P}^{\est}_\Phi(X,\Theta)$ as the ring generated by $P^{\est}_\Phi(X,\Theta)$.
\tria \end{definition}

In summary, we then have the following diagram:
\begin{equation*}
  \begin{tikzcd}
    P_\Phi(X,\theta) \arrow[hookrightarrow]{r} &\tilde{P}_\Phi(X,\theta) \\
    P^{\est}_\Phi(X,\theta) \arrow[hookrightarrow]{r} \arrow[hookrightarrow]{u} & \tilde{P}^{\est}_\Phi(X,\theta). \arrow[hookrightarrow]{u}
  \end{tikzcd}
\end{equation*}

\begin{example}
  Let $(N,\theta)$, $\Phi$, and $f=\frac{1}{x_1x_2}$ be as in Example \ref{potentialeg}. We showed that $f$ is dominated by $\Phi$, but it is in fact estimate-dominated: We have $\left\{\frac{1}{x_1},\frac{x_1+x_3}{x_2x_3}\right\}\subset \Phi$, and
  \[
    \left(\frac{1}{x_1}\right)\left(\frac{x_1+x_3}{x_2x_3}\right)-\frac{1}{x_1x_2}=\frac{1}{x_2x_3}\in P(N,\theta).
  \]
\end{example}

\begin{remark}
  In general, the inclusion $P_\Phi^{\est}(X,\Theta)\hookrightarrow P_\Phi(X,\Theta)$ is strict. For instance, let $(X,\Theta)=(\G_m,[\Id])$ and $\Phi=\{x^2\}$. Then $\Ccal_x(\G_m,[\Id])=\Ccal_\Phi(\G_m,[\Id])$, but there is no positive rational function $g$ on $\G_m$ so that $x+g$ is a polynomial in $x^2$.

  The same example shows that the assignments $P_\bullet^{\est}$ and $\tilde{P}_\bullet^{\est}$ are not functorial. The identity map
  \[
    \Id: (\G_m,[\Id],x^2)\to (\G_m,[\Id],x)
  \]
  is an isomorphism of positive varieties with potential, but we saw $P_x^{\est}(\G_m,[\Id])\ne P_{x^2}^{\est}(\G_m,[\Id])$.
\end{remark}

\begin{remark}
  In \cite{BK}, the author consider a single potential function $\Phi=\{f\}$. For our purposes it is convenient to consider sets of potentials. This is equivalent to considering a single potential, as shown in the following proposition.
\end{remark}

\begin{prop}
\label{singlepotential}
  Let $(X,[\theta])$ be a positive variety with set of potentials $\Phi=\{\varphi_1,\dots, \varphi_m\}$, and let $\Phi^+=\{\varphi_1+\cdots+\varphi_m\}$. Then $P_{\Phi}(X,[\theta])=P_{\Phi^+}(X,[\theta])$, and $P^{\est}_{\Phi}(X,[\theta])=P^{\est}_{\Phi^+}(X,[\theta])$.
\end{prop}

\begin{proof}
  Since $\Ccal_{\{f,f'\}}(X,\theta)=\Ccal_{f+f'}(X,\theta)$ for $f,f'\in P(X,\theta)$, we have $P_{\Phi}(X,\theta)=P_{\Phi^+}(X,\theta)$.
    
  Since it is obvious that $P^{\est}_{\Phi}(X,[\theta])\supseteq P^{\est}_{\Phi^+}(X,[\theta])$, all we need to show is $P^{\est}_{\Phi}(X,[\theta])\subset P^{\est}_{\Phi^+}(X,[\theta])$. For $f\in P^{\est}_{\Phi}(X,[\theta])$, there exists some $ g\in P(X,[\theta])$ such that $f+g=p$, where $p\in \Poly^+ \Phi$. It suffices to show that there is an $h\in P(X,[\theta])$ such that $p(\varphi_1,\dots,\varphi_m)+h=q$ for $q\in\Poly^+\Phi^+$; then $f+g+h=q$. Without loss of generality we may assume that $p=\varphi_1^{n_1}\cdots \varphi_m^{n_m}$ is a non-constant monomial in the $\varphi_i$'s. Set $q:=(\varphi_1+\cdots+ \varphi_m)^{n_1+\cdots+ n_m}$. Then $q-p\in P(X,[\theta]),$	which proves the proposition.
\end{proof}

\section{Potentials on Double Bruhat Cells}

In this Section we recall that double Bruhat cells in semisimple Lie groups are examples of positive varieties. In particular, we focus on the double Bruhat cell $G^{e,w_0}$ and study the semifield of functions weakly dominated by the Berenstein-Kazhdan potentials. Our main technical result is Theorem \ref{maintheorem}, which provides a source of weakly dominated functions to be used in the next Section.

\subsection{Semisimple Groups}\label{semisimple}

Let $G$ be a simply connected semisimple algebraic group over $\Q$ with Lie algebra $\g$ of rank $r$. Choosing a Cartan subalgebra $\hhh\subset \g$ gives a root system $R\subset \hhh^*$; for a root $\alpha\in R$ we denote by $\g_\alpha\subset \g$ the associated root space. Fixing a choice of positive roots $R^+\subset R$ gives a Cartan decomposition $\g=\mathfrak{n}_-\oplus \hhh \oplus \mathfrak{n}$, where 
\[
  \mathfrak{n}_-=\bigoplus_{-\alpha\in R^+} \g_{\alpha}\qquad
  \mathfrak{n}_+=\bigoplus_{\alpha\in R^+} \g_\alpha.
\] 
Let $N_-$, $H$, and $N$ be closed subgroups of $G$ with Lie algebras $\mathfrak{n}_-, \hhh,$ and $\mathfrak{n}$ respectively. Then $H$ is a Cartain subgroup of $G$ and $N_-$, $N$ are a pair of opposite maximal unipotent subgroups of $G$. Let $F_i,H_i,E_i$ be the Chevalley generators of $\g$, indexed by the simple roots $\alpha_1,\dots,\alpha_r\in \hhh^*$ of $\g$. For each $i$, the triple $F_i,H_i,E_i$ defines an embedding $\phi_i: \rm{SL}_2\to G$.
We define
\[
  f_i(t):=\phi_i \begin{pmatrix} 1 & 0 \\ t & 1 \end{pmatrix},
  \quad h_i(t):=\phi_i \begin{pmatrix} t & 0 \\ 0 & -t \end{pmatrix},\quad e_i(t):=\phi_i \begin{pmatrix} 1 & t \\ 0 & 1 \end{pmatrix}.
\]
Let $x\mapsto x^T$ be the involutive \emph{transpose} anti-automorphism of $G$ defined by
\[
  h_i(t)^T=h_i(t),~e_i(t)^T=f_i(t),~f_i(t)^T=e_i(t).
\]

We denote by $W$ the Weyl group $N_G(H)/H$ of $G$. There is a faithful linear action of $W$ on $\hhh^*$: identifying $W$ with its image under this action; the group $W$ is generated by the simple reflections $s_1,\dots,s_r$, where $s_i(\omega)= \omega-\omega(H_i)\alpha_i$.
For the simple reflection $s_i\in W$, set the lift
\[
  \overline{s}_i=\phi_i\begin{pmatrix} 0 & -1 \\ 1 & 0 \end{pmatrix}
\]
of $s_i$ to $G$. It is known that the $\overline{s}_i$'s satisfy the Coxeter relations in $W$, thus any reduced word of $w\in W$ gives the same lift $\overline{w}\in G$.

Consider the opposite Borel subgroups $B=HN$ and $B_-=HN_-$, with Lie algebras $\mathfrak{b}=\hhh\oplus \mathfrak{n}_+$ and $\mathfrak{b}_-=\hhh\oplus \mathfrak{n}_-$, respectively. Then $G$ has a decomposition into \emph{double Bruhat cells}:
\[
  G=\bigsqcup_{u,v} G^{u,v}, \qquad \text{where~} G^{u,v}:=B u B \cap B_- v B_-.
\]
The expression on the right does not depend on the lift of $u,v\in W$ to $N_G(H)\subset G$.

\begin{example}
\label{sl3liegroupexample} We introduce an example which we will reference in future sections.
    Let $G={\rm SL}_3$, then the Lie algebra of $G$ is $\mathfrak{sl}_3 \subset \Mat (3\times 3)$, traceless $3\times 3$-matrices. Let $E_{ij}\in \Mat (3\times 3)$ denote the matrix with all entries zero except the $ij$-th entry:
\[ 
(E_{ij})_{kl}=\delta_{ik}\delta_{jl}.
\]
For $i\in\{1,2\}$, the simple coroots are $H_i=E_{ii}-E_{i+1,i+1}$. We make the standard choice of $E_{12}$ and $E_{23}$ as simple root vectors of $\mathfrak{sl}_3$. Then $B\subset {\rm SL}_3$ is the subgroup of upper triangular matrices and $B_-\subset {\rm SL}_3$ is the subgroup of lower triangular matrices. In this case, the transpose map $x\mapsto x^T$ defined above is the usual matrix transpose. The Weyl group of ${\rm SL}_3$ is isomorphic to the symmetric group $S_3$, with simple reflections corresponding to the transpositions $s_1=(12)$ and $s_2=(23)$, respectively. 
\end{example}

The \emph{weight lattice} $P(\g)$ of $\g$ is the set of $\omega\in \hhh^*$ which have $\omega(H_i)\in\Z$ for all $i$. There is a basis of \emph{fundamental weights} $\omega_i$ of $P$, where $\omega_i(H_j)=\delta_{ij}$. Any multiplicative character $\xi:H\to \G_m$ induces a weight $\omega=d\xi\in \mathfrak{h}^*$; by assumption $G$ is simply-connected and so we may identify multiplicative characters $\Hom(H,\G_m)$ with the weight lattice $P(\g)$. We write $h^\omega$ for the value at $h\in H$ of the character corresponding to $\omega\in P(\g)$. For the fundamental weight $\omega_i$ we have $h_j(t)^{\omega_i}=t^{\delta_{ij}}$.

Let $G_0=N_- H N\subset G$ be the open subset of $G$ of elements which admit a Gaussian decomposition. For an element $x\in G_0$ we write $x=[x]_- [x]_0 [x]_+$. 
Following \cite{BZ}, for $u,v\in W$ and a fundamental weight $\omega_i$ we define the \emph{generalized minor} $\Delta_{u\omega_i,v\omega_i}$ as  the regular function on $G$ whose restriction to the open subset $\overline{u} G_0 \overline{v}\n$ of $G$ is given by 
\[
\Delta_{u\omega_i,v\omega_i}(x)=[\overline{u}\n x \overline{v}]_0^{\omega_i}.
\]
It is shown in \cite{FZ} that the right hand side depends only on the weights $u\omega_i$ and $v\omega_i$, and not on the choice of lifts $\overline{u}$ and $\overline{v}$. If $a_1,a_2\in H$ and $x\in G$, then
\begin{equation}
\label{BonusIdentity}
\Delta_{u\omega_i,v\omega_i}(a_1 x a_2)= a_1^{u\omega_i} a_2^{v\omega_i} \Delta_{u\omega_i,v\omega_i}.
\end{equation}
If $G={\rm SL}_n$, the generalized minors are minors. 

\subsection{Positive Structures on Double Bruhat Cells}
 
Let $G$ be as in the previous section, and let $W$ be its Weyl group. Let $(u,v)\in W\times W$. Recall that $W$ is generated by the simple reflections $s_i$ corresponding to the simple roots of $\g$, and that a word $u$ in the simple reflections of $W$ has a length $l(u)$. We write $s_{-1},\dots,s_{-r}$ for the simple reflections in the first copy of $W$ and $s_1,\dots,s_r$ for the simple reflections in the second copy of $W$. A \emph{double word} $\mathbf{i}$ for $(u,v)$ is a shuffle of a word for $u$ in the simple reflections $s_{-1},\dots,s_{-r}$, and a word for $v$ in the simple reflections $s_1,\dots,s_r$. If the words for $u$ and $v$ are both reduced, then $\mathbf{i}$ is a \emph{double reduced word} for $(u,v)$. The \emph{length} of $\mathbf{i}$ is  $l(\mathbf{i}):=l(u)+l(v)$. We write
\[
  \mathbf{i}=(i_1,\dots,i_{l(u)+l(v)}),
\]
where $i_j\in\{-r,\dots,-1,1,\dots,r\}$. 
 
A double word defines a map $\theta_{\mathbf{i}}:  \G_m^{l(\mathbf{i})+r}\to G$. The map is given by:
\begin{equation}\label{rep}
  (t_1,\dots,t_{l(\mathbf{i})},t_{l(\mathbf{i})+1},\dots, t_{l(\mathbf{i})+r})\mapsto e_{i_1}(t_1)\cdots e_{i_{l(\mathbf{i})}}(t_{l(\mathbf{i})})\cdot h_1(t_{l(\mathbf{i})+1})\cdots h_r(t_{l(\mathbf{i})+r})
\end{equation}
where we interpret $e_{-i}(t):=f_i(t)$, for $i>0$.

\begin{example}
\label{sl3factparexample}
  Let $G={\rm SL}_3$ as in Example \ref{sl3liegroupexample}, and fix $(u,v)=(e,w_0)$, where $w_0=(13)$ is the longest element of $S_3$, the Weyl group of $G$. A double reduced word for $(e,w_0)$ is $\mathbf{i}=s_1 s_2 s_1$. Then $\theta_{\mathbf{i}}:\G_m^3\times \G_m^2 \to {\rm SL}_3$ has
  \begin{align*}
    \theta_{\mathbf{i}}(t_1,t_2,t_3,t_4,t_5) &= \theta(t_1,t_2,t_3)\cdot h_1(t_4)\cdot h_2(t_5),
  \end{align*}
  where $\theta:\G_m^3\to {\rm SL}_3$ is given as in Example \ref{charteg1}.
\end{example}

\begin{theorem}[Theorem 1.2 of \cite{FZ}]
\label{FZthm1.2}
  For any $(u,v)\in W\times W$ and any double reduced word $\mathbf{i}$ for $(u,v)$ the map $\theta_{\mathbf{i}}$ restricts to a biregular isomorphism between $\G_m^{l(\mathbf{i})+r}$ and a Zariski open subset of $G^{u,v}$.
\end{theorem}

As a consequence of the previous theorem and Theorem 3.1 from \cite{BZ}, we have the following:

\begin{theorem}\label{isoposstruc}
  Let $\mathbf{i}$ and $\mathbf{i'}$ be double reduced words for $(u,v)\in W\times W$. Then $\theta_{\mathbf{i}}$ and $\theta_{\mathbf{i'}}$ are positively equivalent toric charts on $G^{u,v}$.
\end{theorem}

For the remainder of this article, we restrict our attention to the case of $(u,v)=(e,w_0)$, where $w_0$ is the longest element of $W$. A double reduced word for $(e,w_0)$ is nothing but a reduced word for $w_0$. We have the following relationship between generalized minors and the toric charts $\theta_\mathbf{i}$. From  Theorem 5.8(i) and formula (4.2), both of \cite{BZ2}, we have

\begin{theorem}[]\label{thm5.8}
	 Let $\gamma$ and $\delta$ be two weights in the $W$-orbit of the same fundamental weight $\omega_i$ of $\mathfrak{g}$, and let
	 $\mathbf{i}=(i_1,\dots,i_m)$ be any word in $\{1,\dots,r\}$ and define $\theta_\mathbf{i}:\G_m^{m+r}\to G$ as in (\ref{rep}). Then,
	 \[
	   \Delta_{\gamma,\delta}(\theta_{\mathbf{i}}(t_1,\dots,t_{m+r}))= \Delta_{\gamma,\delta}(\theta_{\mathbf{i}}(t_1,\dots,t_m,1,\dots,1)) \cdot (h_1(t_{m+1})\cdots h_r(t_{m+r}))^\delta,
	 \]
	 where the first term is a polynomial in $t_1,\dots,t_{m}$ with positive integer coefficients and the second term is a Laurent monomial in $t_{m+1}, \dots, t_{m+r}$. In particular, by Theorem \ref{FZthm1.2}, we get $\Delta_{\gamma,\delta}\in P(G^{e,w_0},[\theta_\mathbf{i}])$.
 \end{theorem}

\subsection{Cluster Variables on Double Bruhat Cells}\label{clustersection}

Given a reduced word $\mathbf{i}$ for $w_0 \in W$, there is a set $F(\mathbf{i})\subset P(G^{e,w_0},\theta_\mathbf{i})$ of generalized minors
\[
  \Delta(-r;\mathbf{i}),\dots,\Delta(-1;\mathbf{i}),\Delta(1;\mathbf{i}),\dots,\Delta(l(w_0);\mathbf{i}),
\]
which are the \emph{cluster variables} associated to the \emph{initial seed} of $\mathbb{Q}[G^{e,w_0}]$ determined by $\mathbf{i}$; details are given in \cite{BFZ}. 

\begin{example}
\label{sl3posstrucexample}
  Let $G={\rm SL}_3$ as in Example \ref{sl3liegroupexample} and \ref{sl3factparexample}.
  Consider the reduced word $\mathbf{i}=s_1s_2s_1$ for $w_0$, where $w_0=(13)$ is the longest element of of $S_3$. Then the set $F(\mathbf{i})$ of cluster variables for ${\rm SL}_3^{e,w_0}$ is
\[ 
F(\mathbf{i})=\{ 
\Delta_{12,23},\Delta_{1,3},\Delta_{1,2},\Delta_{12,12},\Delta_{1,1}
\}, 
\] 
where $\Delta_{I,J}\in\Q[\rm{SL}_3]$ is the minor with rows $I$ and columns $J$.
\end{example}

\begin{remark}\label{frozenvariables}
  For any choice of reduced word $\mathbf{i}$ for $w_0$, we always have $\Delta(-k,\mathbf{i})=\Delta_{\omega_k,w_0\omega_k}\in F(\mathbf{i})$. Also, we always have the \emph{principal} generalized minors $\Delta_{\omega_k,\omega_k}\in F(\mathbf{i})$, for all $k$. Together, these are the \emph{frozen variables} of the upper cluster algebra structure on $\Q[G^{e,w_0}]$.
\end{remark}

Define the rational map
\begin{align*}
    \Delta_{\mathbf{i}}: G^{e,w_0} & \to \G_m^{r+l(w_0)}, \\
    g&\mapsto (\Delta(-r;\mathbf{i})(g),\dots,\Delta(-1;\mathbf{i})(g),\Delta(1;\mathbf{i})(g),\dots,\Delta(l(w_0);\mathbf{i})(g)).
\end{align*}
\begin{prop}[Lemma 2.12 of \cite{BFZ}]
     \label{clusterposstructure}
     The map $\Delta_\mathbf{i}$ is biregular isomorphism from Zariski open subset of $G^{e,w_0}$ to $\G_m^{r+l(w_0)}$.
\end{prop}
We then have, for each reduced word $\mathbf{i}$ for $w_0$, a toric chart
\[
  \Delta_\mathbf{i}\n:\G_m^{r+l(\mathbf{i})}\to G^{e,w_0}.
\]

\begin{prop}\label{twistmap} For any reduced word $\mathbf{i}$ for $w_0$, the maps $\Delta_\mathbf{i}\n$ and $\theta_\mathbf{i}$ are positively equivalent toric charts on $G^{e,w_0}$.
\end{prop}

\begin{proof}
  By Theorem \ref{thm5.8}, the change of coordinates $\Delta_\mathbf{i}\circ \theta_\mathbf{i}$ is positive and birational. We need to show its inverse $\theta\n_\mathbf{i}\circ \Delta_\mathbf{i}\n$ is positive. In \cite{FZ}, the authors introduce an involutive biregular \emph{twist map} $\xi:G^{e,w_0}\to G^{e,w_0}$, and Theorem 1.9 in \cite{FZ} gives that $\theta_\mathbf{i}$ and $\xi\circ \Delta_\mathbf{i}\n$ are positively equivalent toric charts on $G^{e,w_0}$. So,
  \[
    \theta_\mathbf{i}\n\circ \Delta_\mathbf{i}\n=(\theta_\mathbf{i}\n\circ \xi\circ \Delta_\mathbf{i}\n)\circ(\Delta_\mathbf{i}\circ \theta_\mathbf{i})\circ(\theta_\mathbf{i}\n\circ \xi\circ \Delta_\mathbf{i}\n)
  \]
  is positive.
\end{proof}

In summary, for any two reduced words $\mathbf{i}, \mathbf{i}'$ for $w_0$, the following positive structures on $G^{e,w_0}$ are all equal:
\begin{equation}\label{equivequation}
  [\theta_{\mathbf{i}}]=[\theta_{\mathbf{i}'}]=[\Delta\n_\mathbf{i}]=[\Delta\n_{\mathbf{i}'}].
\end{equation}

\subsection{Weakly Estimate-Dominated Functions on Double Bruhat Cells} \label{wedfodbc}

Let $\mathbf{i}$ be a reduced word for $w_0$. In this section we consider $G^{e,w_0}$ with the positive structure determined by $\theta_\mathbf{i}$. 

Recall that the Lie algebra $\g$ of $G$ has a left action on $\Q[G]$; in particular for the elements $E_i,F_i$ of the Chevalley basis of $\g$, we have for $\Delta\in \Q[G]$,
\[
  F_i\cdot \Delta(g)=\frac{d}{dt}\Big|_{t=0} \Delta(f_i(t)g), \quad
  E_i\cdot \Delta(g)=\frac{d}{dt}\Big|_{t=0} \Delta(e_i(t)g).
\]
If $X\in \nnn$, then $X$ is everywhere tangent to the double Bruhat cell $G^{e,w_0}$. Define the potentials
\[
  \varphi_i:= \frac{E_i\cdot\Delta_{\omega_i,w_0\omega_i}}{\Delta_{\omega_i,w_0\omega_i}}=\frac{\Delta_{s_i\omega_i,w_0\omega_i}}{\Delta_{\omega_i,w_0\omega_i}},
\]
and let $\Phi:=\{\varphi_1,\dots,\varphi_{r}\}$. The following proposition implies that $\Phi\subset P(G^{e, w_0},[\theta_\mathbf{i}])$.

\begin{prop}[Proposition 4.11.(ii) in \cite{BZ2}]\label{pro4.11}
 	Let $\mathbf{i}=(i_1,\dots,i_m)$ be a reduced 
 	word for $w_0$. Then 
 	\[
 		\frac{E_{i_1}\cdot \Delta_{\omega_{i_1},w_0\omega_{i_1}}(\theta_\mathbf{i}(t_1,\dots,t_m,t_{m+1},\dots,t_{m+r}))}{\Delta_{\omega_{i_1},w_0\omega_{i_1}}( \theta_\mathbf{i}(t_1,\dots,t_m,t_{m+1},\dots,t_{m+r}))}=\frac{1}{t_1}.
 	\]
\end{prop}

Note that on the toric chart $\Delta_\mathbf{i}\n$, the potentials $\varphi_j$ restrict to regular functions and so, by remark \ref{coneremark}, the cone $\Ccal_\Phi(G^{e,w_0},\Delta_\mathbf{i}\n)(\R)$ is convex. Since $\Q[G]$ is a $\mathfrak{g}$ module, it is a $U_{\mathfrak{g}}$ module. Next we will estimate the action of $U_{\mathfrak{g}}$. In fact, we have 

\begin{lemma}\label{mainlemma}
  For a word $\mathbf{j}=(j_1,\dots,j_n)$ in the simple roots of $\g$, choose a reduced word $\mathbf{i}=(i_1,\dots,i_m)$ for $w_0$ such that $i_1=j_1$. Suppose that $E_{j_2}\cdots E_{j_n}\Delta_{\omega_{i},w\omega_{i}}\neq0$. Then for any $w$ in the Weyl group $W$ and any fundamental weight $\omega_i$,
  \[
    \frac{E_{j_1}\cdots E_{j_n}\Delta_{\omega_{i},w\omega_{i}}}{E_{j_2}\cdots E_{j_n}\Delta_{\omega_{i},w\omega_{i}}}\in P^{\est}_\Phi(G^{e,w_0},[\theta_{\mathbf{i}}]).
  \] 
\end{lemma}

\begin{proof}
	For a word $\mathbf{j}=(j_1,\dots,j_n)$, let $\partial_\mathbf{j}:=\frac{d}{dq_1}\Big|_0\cdots \frac{d}{dq_n}\Big|_0$ be the differential at zero. In the following, let $\mathbf{j}':=(j_2,\dots, j_n)$ and $\theta_\mathbf{i}$ be $\theta_\mathbf{i}(t_1,\dots,t_m,t_{m+1},\dots,t_{m+r})$ as in \eqref{rep}. By Theorem \ref{thm5.8}, we get
	\begin{align*}
		E_{j_2}\cdots E_{j_n}\Delta_{\omega_{i},w\omega_{i}}(\theta_\mathbf{i})&=\partial_{\mathbf{j}'} \Delta_{\omega_{i},w\omega_{i}}(e_{j_n}(q_n)\cdots e_{j_2}(q_2)\theta_\mathbf{i}(t_1,\dots,t_{m+r}))\\
		&=\partial_{\mathbf{j}'} \Delta_{\omega_{i},w\omega_{i}}(e_{j_n}(q_n)\cdots e_{j_2}(q_2)e_{i_1}(t_1)\cdots e_{i_m}(t_m)h_1(t_{m+1})\cdots h_r(t_{m+r}))\\
		&=\sum_{k=0} f_k(t_2,\dots, t_{m+r})t_1^k,
	\end{align*}
	where $f_k$ are products of polynomials in $t_2, \dots, t_m$ with positive integer coefficients and Laurent monomials in $t_{m+1}, \dots, t_{m+r}$. Similarly, since $j_1=i_1$,
	\begin{align*}
		E_{j_1}\cdots E_{j_n}\Delta_{\omega_{i},w\omega_{i}}(\theta_\mathbf{i})&=\partial_\mathbf{j} \Delta_{\omega_{i},w\omega_{i}}(e_{j_n}(q_{n})\cdots e_{j_1}(q_{1})\theta_\mathbf{i}(t_1,\dots,t_{m+r}))\\
		&=\partial_\mathbf{j} \Delta_{\omega_{i},w\omega_{i}}\big(e_{j_n}(q_{n})\cdots e_{j_2}(q_2) e_{i_1}(q_{1}+t_{1})e_{i_2}(t_{2})\cdots \\
		& \qquad \qquad \cdots e_{i_m}(t_m)h_1(t_{m+1})\cdots h_r(t_{m+r})\big)\\
		&=\frac{d}{dq_{1}}\Big|_0\sum_{k=0} f_k(t_{2},\dots, t_{m+r})(q_{1}+t_{1})^k\\
		&=\sum_{k=1} k f_k(t_{2},\dots, t_{m+r})t_{1}^{k-1}.
	\end{align*}
	By Proposition \ref{pro4.11}, we know $\varphi_{i_
	1}(\theta_\mathbf{i}(t_1,\dots,t_{m+r}))=1/t_{1}$. Direct calculation gives:
	\[
		\frac{\alpha}{t_{1}}-\frac{E_{j_1}\cdots E_{j_n}\Delta_{\omega_{i},w_0\omega_{i}}(\theta_\mathbf{i}(t_1,\dots,t_{m+r}))}{E_{j_2}\cdots E_{j_n}\Delta_{\omega_{i},w_0\omega_{i}}(\theta_\mathbf{i}(t_1,\dots,t_{m+r}))}=\frac{\alpha f_0+\sum_{k=1} (\alpha-k)f_kt_{1}^{k}}{t_1\sum_{k=0} f_kt_{1}^k}. 	
	\]
	Let $\alpha$ be a sufficiently large positive integer, then the right hand side is a positive function. 
\end{proof}

\begin{remark}   \label{totally_positive}
In fact, the arguments in the proof of Lemma \ref{mainlemma} are applicable in a much larger context. We say that a function $f$  on $G^{e, w_0}$ is {\em totally positive} if for every word ${\bf j} = (j_1, \dots, j_n)$ the expression $E_{j_1} \dots E_{j_n} f$ is either a positive rational function or zero. The functions $\Delta_{\omega_i, w_0 \omega_i}$ satisfy this property, but there are many other examples: for $G$ simply laced, all elements of the dual canonical basis have this property, and for $G$ semisimple all cluster variables on $G^{e, w_0}$ are totally positive. 

If $f$ is totally positive, $E_{j_2} \dots E_{j_n} f$ (and all such derivatives) are products of polynomials with positive coefficients in  variables $t_1, \dots, t_m$ and of Laurent monomials in $t_{m+1}, \dots, t_{m+r}$. Hence, the proof of Lemma \ref{mainlemma} applies and the logarithmic derivative 
$$
\frac{E_{j_1}E_{j_2} \dots E_{j_n} f}{E_{j_2} \dots E_{j_n} f}
$$
is dominated by the potential $\phi_{j_1}=1/t_1$. 

The theory of totally positive functions and its applications to functions dominated by potentials will be explored elsewhere.
\end{remark}

\begin{theorem}\label{maintheorem}
	For any nonzero vector $X\in \g_\alpha$ in a positive root space of $\g$ and any $w\in W$, we have $\frac{X\cdot \Delta_{\omega_i,w\omega_i}}{\Delta_{\omega_i,w\omega_i}}\in\tilde{P}_\Phi^{est}(G^{e,w_0},[\theta_\mathbf{i}])$.
\end{theorem}

\begin{proof}
  Without loss of generality, we represent $X$ as a nested commutator of simple roots $X=[E_{k_1}[\cdots[E_{k_{n-1}},E_{k_n}]\cdots ]]$. Since $\tilde{P}_\Phi^{est}(G^{e,w_0},[\theta_\mathbf{i}])$ is a ring, it suffices to show that for any word $\mathbf{j}=(j_1,\dots,j_n)$ in the simple roots of $\g$,
  \[
  \frac{E_{j_1}\cdot E_{j_2}\cdots E_{j_n} \Delta_{\omega_i,w\omega_i}}{E_{j_2}\cdots E_{j_n}  \Delta_{\omega_i,w\omega_i}}\in P_\Phi^{est}(G^{e,w_0},[\theta_\mathbf{i}]).
  \]
  Since any choice of a reduced word $\mathbf{i'}$ for $w_0$ defines  toric chart $\theta_{\mathbf{i}'}$ on $G^{e,w_0}$ which is positively equivalent to $\theta_{\mathbf{i}}$, this follows immediately from Lemma \ref{mainlemma}.
\end{proof}


An analogous statement holds for the standard right action of $\g$ on $\C[G]$. For the simple root $\alpha_i'=-w_0\alpha_i$, write $E_{i'}$ for the corresponding root vector. Define the potentials $\psi_i:= \frac{\Delta_{\omega_i,w_0\omega_i}\cdot E_{i'}}{\Delta_{\omega_i,w_0\omega_i}}=\frac{\Delta_{\omega_i,w_0s_i\omega_i}}{\Delta_{\omega_i,w_0\omega_i}}$, and let $\Psi:=\{\psi_1,\dots,\psi_{r}\}$. Then $\Psi\subset P(X,[\theta_\mathbf{i}])$, and

\begin{theorem}\label{maintheoremprime}
  For any nonzero vector $X\in \g_\alpha$ in a positive root space of $\g$ and any $w\in W$, we have $\frac{\Delta_{\omega_i,w\omega_i}\cdot X}{\Delta_{\omega_i,w\omega_i}}\in\tilde{P}_\Psi^{est}(G^{e,w_0},[\theta_\mathbf{i}])$.
\end{theorem}

The proof is symmetric to the proof of Theorem \ref{maintheorem}.

\begin{definition}
\label{BKpotential}
  For any reduced word $\mathbf{i}$ for $w_0$, consider the positive variety with potential $(G^{e,w_0},[\theta_\mathbf{i}],\Phi_{BK})$, where 
  \[ \Phi_{BK}:=\{ \varphi_1,\dots,\varphi_r,\psi_1,\dots,\psi_r\}
  \]
  with $\varphi_i,\psi_i$ are as above. We call $\Phi_{BK}$ the \emph{Berenstein-Kazhdan Potential} on $G^{e,w_0}$. For a certain toric chart $\theta_z \in[\theta_\mathbf{i}]$ on $G^{e,w_0}$, the cone $\Ccal^\le_{\Phi_{BK}}\subset(G^{e,w_0},\theta_z)^t$ is called the \emph{extended string cone}. Abusing terminology slightly, we call $\Ccal_{\Phi_{BK}}\subset(G^{e,w_0},\theta)^t$ the \emph{strict extended string cone} for any choice of $\theta\in[\theta_\mathbf{i}]$; any such cone is related to $\Ccal_{\Phi_{BK}}\subset (G^{e,w_0},\theta_z)^t$ by a piecewise linear bijection.
  See \cite{BK} for more details, including the connection between $\Ccal^\le_{\Phi_{BK}}$ and parameterizations of crystal bases. The potential $\Phi_{BK}$ was also discovered independently by Rietsch, see \cite{Ri}. For $G={\rm SL}_n$, for a certain choice of reduced word $\mathbf{i}$ for $w_0$, the cone $\Ccal^\le_{\Phi_{BK}}\subset(G^{e,w_0},\theta_\mathbf{i})^t$ is the Gelfand-Zeitlin cone \cite{BZ1}.
\tria \end{definition}

\section{Positive Poisson Varieties} 
\label{lastsection}

In this Section we consider Poisson structures on positive varieties. We are interested in weakly log-canonical charts. These are toric charts in which the Poisson structure is log-canonical modulo terms weakly dominated by potentials. Examples are given by dual Poisson algebraic groups $G^*$ for $G$ a semisimple algebraic group.

\subsection{Definition of Positive Poisson Varieties}

Let $(X,\pi)$ be an irreducible Poisson variety over $\Q$, and let $\Theta$ be a positive structure on $X$. Suppose $\theta:\G_m^n\to X$ is a toric chart with $\theta\in \Theta$ and let $(z_1,\dots,z_n)$ be the standard coordinates on $\G_m^n$. Fix a set of potentials $\Phi\subset P(X,[\theta])$ on $X$. Then $\theta$ is a \emph{weakly log-canonical chart} for $\pi$ if, in the coordinates $z_i$ the bracket is of the form 
\begin{equation}\label{logcanon}
  \{z_i,z_j\}=z_iz_j(\pi_{ij}+ f_{ij}),
\end{equation}
where $\pi_{ij}\in k$ are constant and $f_{ij}\in \tilde{P}_{\Phi}(X, [\theta])$ are weakly dominated by $\Phi$. We call $\pi_{ij}$ the \emph{log-canonical part} of the bracket.

If a Poisson variety $(X,\pi)$ has a weakly log-canonical chart $\theta$, we define $[\theta]_\pi$ to be the collection of all $\theta'\in [\theta]$ that is weakly log-canonical chart for $\pi$.  We call $[\theta]_\pi$ a $\pi$-\emph{compatible positive structure} on $X$.


We define a \emph{positive Poisson variety} to be a quadruple $(X,\pi,\Theta_{\pi},\Phi)$, where $\Theta_\pi$ is a $\pi$-compatible positive structure on $X$. Note that $[\theta]_\pi\subset [\theta]$ but in general $[\theta]\ne[\theta]_\pi$.

A \emph{positive Poisson map} of positive Poisson varieties
\[
  \phi: (X_1,\pi_1,\Theta_{\pi_1},\Phi_1)\to (X_2,\pi_2,\Theta_{\pi_2},\Phi_2)
\]
is a Poisson map $\phi: (X_1,\pi_1)\to (X_2,\pi_2)$, which is also a map of positive varieties with potential. We denote by $\mathbf{PosPoiss}$ the category of positive Poisson varieties.

\begin{example}
    In \cite{GSV}, the authors define the \emph{cluster manifold} $X(A)$ associated to a cluster algebra $A$. The cluster algebra structure on $A$ gives $X(A)$ a Poisson structure $\pi$ and a family of positively equivalent toric charts $\theta$ on $X$ which are log-canonical for $\pi$. Cluster manifolds are then examples of positive Poisson varieties, with potential $\Phi=0$; see Remark \ref{coneremark}.
\end{example}

We record the following observation for future reference.

\begin{prop}\label{monomiallogcan}
  Let $\theta:\G_m^n\to X$ be a weakly log-canonical chart for $(X,\pi,\Theta_\pi,\Phi)$. If $M$ and $N$ are Laurent monomials in the standard coordinates $z_j$, then the bracket $\{M,N\}$ is weakly log-canonical; i.e.
  \[
    \{M,N\}=MN(\pi_{MN}+f_{MN}),
  \]
  where $\pi_{MN}\in k$ and $f_{MN}\in \tilde{P}_\Phi(X,[\theta])$.
\end{prop}

\begin{proof}
  Assume $M,N$, and $L$ are Laurent monomials in $z_1,\dots,z_n$ and that
  \[
    \{M,N\}=MN(\pi_{MN}+f_{MN}) \mbox{~and~} \{M,L\}=ML(\pi_{ML}+f_{ML})
  \]
  are weakly log-canonical. The proposition follows by induction, using the following two facts. First,
  \[
    \{M,N\n\}=-N^{-2}\{M,N\}=N^{-2}MN(-\pi_{MN}-f_{MN})=MN\n(-\pi_{MN}-f_{MN})
  \]
  is weakly log-canonical. Second,
  \begin{align*}
    \{M,NL\} & =\{M,N\}L+N\{M,L\} = MNL(\pi_{MN}+f_{MN})+MNL(\pi_{ML}+f_{ML})
  \end{align*}
  is weakly log-canonical.
\end{proof}

\subsection{Poisson Algebraic Groups}\label{PLG}

Most of the material in this section is based on \cite{ES,LW}. Let $G$ be a simply connected semisimple algebraic group over $\Q$ with Lie algebra $\g$, as in Section \ref{semisimple}. Let $X_i$ be an orthonormal basis for $\mathfrak{h}$ under the Killing form. For positive roots $\alpha\in R^+$, choose root vectors $E_\alpha\in \g_\alpha$ and $E_{-\alpha}\in \g_{-\alpha}$ so that under the Killing form $K(\cdot,\cdot)$ we have $K(E_\alpha,E_{-\alpha})=1$. Consider the standard quasitriangular $r$-matrix
\[
  r_G:=\frac{1}{2}\sum_{i=1}^{\operatorname{rank} \hhh} X_i\otimes X_i +
    \sum_{\alpha\in R^+}E_{\alpha}\otimes E_{-\alpha}.
\]
Let $r_G^\lambda,r_G^\rho$ be the left- and right-invariant 2-tensor fields on $G$, respectively, which have $r_G^\lambda(e)=r_G^\rho(e)=r_G$. Then $r_G$ satisfies the classical Yang-Baxter equation, and so $\pi_G=r_G^\lambda-r_G^\rho$ is a Poisson bivector field and $(G,\pi_G)$ is a Poisson algebraic group. The simply connected Drinfeld double of $(G,\pi_G)$ is $(D,\pi_D)=(G\times G,\pi_D)$, where $\pi_D=r_D^\lambda-r_D^\rho$ and
\[
  r_D=\frac{1}{2} \sum_{i=1}^{\operatorname{rank} \hhh} (X_i,X_i)\otimes (X_i,-X_i) +\sum_{\alpha\in R^+} (E_\alpha,E_\alpha)\otimes (0,-E_{-\alpha}) +(E_{-\alpha},E_{-\alpha})\otimes (E_\alpha,0)
\]
is the standard quasitriangular $r$-matrix. The diagonal embedding $(G,\pi_G)\to(D,\pi_D)$, sending $g\mapsto(g,g)$ is Poisson.

Recall that $[\cdot]_0:G_0\to H$ is the projection of Gaussian decomposable elements to $H$. Then the subgroup $G^*$ of $G\times G$ given by 
\[
   G^*=\{(b^+,b^-)\in B\times B_-| [b^+]_0[b^-]_0=1\}\subset G\times G
\]
is a Poisson algebraic subgroup of $(D,\pi_D)$, and $(G^*,-\pi_D|_{G^*})$ is a Poisson group dual to $(G,\pi_G)$.


\begin{prop}\label{borelpoisson}
  The Borel subgroups $B$ and $B_-$ are Poisson algebraic subgroups of $G$, and the projections
  \[
    \pr_1:G^*\to B,\qquad \pr_2: G^* \to B_-
  \]
  are anti-Poisson. 
\end{prop}

\begin{proof}
    Follows from the expression of the $r$-matrices $r_G$ and $r_D$; see also \cite{KZ}.
\end{proof}

\begin{definition}[Notation]
\label{projnotation}
  Let $f\in \Q[B]$ and $g\in \Q[B_-]$ be regular functions. Denote $f_1:=f\circ \pr_1\in \Q[G^*]$ and $g_2:=g\circ \pr_2\in \Q[G^*]$.
\tria \end{definition}

For future reference, we note the following computation, which follows directly from the description of $\pi_{G^*}$ above.

\begin{prop}
  Let $f\in \Q[B]$ and $g\in \Q[B_-]$.  Then in the notation of Definition \ref{projnotation},
  \begin{align}\label{mixedbracket}
    \{f_1,g_2\}_{G^*} =\quad &\frac{1}{2}\sum_{i=1}^{\operatorname{rank} \hhh} (X_i\cdot f)_1(-X_i\cdot g)_2 - (f\cdot X_i)_1(g\cdot (-X_i))_2 \\
    +&\sum_{\alpha\in R^+} (E_\alpha\cdot f)_1( -E_{-\alpha}\cdot g)_2 -(f\cdot E_\alpha)_1(g\cdot (-E_{-\alpha}))_2. \nonumber
  \end{align}
\end{prop}

\subsection{The Positive Poisson Variety \texorpdfstring{$G^*$}{G}}

In this section we endow $G^*$ with the structure of a positive Poisson variety. 
Recall the transpose $(\cdot)^T:G\to G$ anti-involution defined in Section \ref{semisimple}. We introduce the Lie group involution $\tau:G\to G$, with $\tau(g)=(g^T)\n=(g\n)^T$. From the expression for the $r$-matrix $r_G$, we have the following.

\begin{prop}\label{tauAP}
  The map $\tau:(G,\pi_G)\to (G,\pi_G)$ is anti-Poisson. Also, $\tau$ restricts to the Lie group isomorphisms $\tau:B\to B_-$ and $\tau:B_-\to B$, and induces the automorphism
  \[
    \tau:G^*\to G^*\ :\  (b^+,b^-)\mapsto (((b^-)^T)\n,((b^+)^T)\n).
  \]
  \end{prop}

Recall from section \ref{clustersection}, we have for each reduced word $\mathbf{i}$ for $w_0$ a set of $r+l(w_0)$ cluster variables $F(\mathbf{i})$ which includes the generalized minors $\Delta_{\omega_k,\omega_k}$. For each $\mathbf{i}$ we defined a birational map $\Delta_\mathbf{i}:G^{e,w_0}\to \G_m^{r+l(w_0)}$, which had as its components the cluster variables $F(\mathbf{i})$. Let $\eta:\G_m^{r+l(w_0)}\to \G_m^{l(w_0)}$ be the projection parallel to the coordinates $\Delta_{\omega_k,\omega_k}$ of $\Delta_\mathbf{i}$. Let
\begin{equation} \label{Deltahat}
  \hat{\Delta}_\mathbf{i}:=\eta\circ \Delta_\mathbf{i}:G^{e,w_0}\to \G_m^{l(w_0)},
\end{equation}
and let $\Delta_0:G^{e,w_0}\to \G_m^r$ be given by
\begin{equation} \label{Deltazero}
  \Delta_0(g)=(\Delta_{\omega_1,\omega_1}(g),\dots,\Delta_{\omega_r,\omega_r}(g)).
\end{equation}

Recall that $G^{e,w_0}\subset B$ is a Zariski open subset. Consider the rational map
\begin{align}
  \bbDelta_\mathbf{i}: G^* & \to \G_m^{l(w_0)}\times \G_m^{l(w_0)}\times\G_m^r \label{bbDelta}; \\
  (b^+,b^-) &\mapsto (\hat{\Delta}_\mathbf{i}(b^+), \hat{\Delta}_\mathbf{i}(\tau(b^-)),\Delta_0(b^+)). \nonumber 
\end{align}

The following proposition follows from the definition of $G^*$ and Proposition \ref{clusterposstructure}.

\begin{prop}
  The map $\bbDelta_\mathbf{i}$ is a biregular isomorphism from a Zariski open subset of $G^*$ to $\G_m^{r+2l(w_0)}$.
\end{prop}

\begin{example}
Let $G={\rm SL}_3$ as in Examples \ref{sl3liegroupexample}, \ref{sl3factparexample}, and \ref{sl3posstrucexample} above. For $b\in G^{e,w_0}$, we have
\begin{align*}
\Delta_{\mathbf{i}}(b) & = (\Delta_{12,23}(b),\Delta_{1,3}(b),\Delta_{1,2}(b),\Delta_{12,12}(b),\Delta_{1,1}(b)), \\ 
\hat{\Delta}_{\mathbf{i}}(b) & = (\Delta_{12,23}(b),\Delta_{1,3}(b),\Delta_{1,2}(b)), \\ 
\Delta_0(b) & = (\Delta_{1,1}(b),\Delta_{12,12}(b)).
\end{align*}
 The birational map $\bbDelta_\mathbf{i}:G^*\to \G_m^8$ is then given:
\begin{align*}
    \bbDelta_\mathbf{i}(b^+,b^-) = & (\Delta_{12,23}(b^+),\Delta_{1,3}(b^+),\Delta_{1,2}(b^+), \\ & \Delta_{12,23}(\tau(b^-)),\Delta_{1,3}(\tau(b^-)),\Delta_{1,2}(\tau(b^-)),\Delta_{1,1}(b^+),\Delta_{12,12}(b^+)).
\end{align*}
\end{example}

Consider the toric chart $\bbDelta_\mathbf{i}\n:\G_m^{r+2l(w_0)}\to G^*$ on $G^*$. By (\ref{equivequation}), another reduced word $\mathbf{i}'$ for $w_0$ gives a toric chart $\bbDelta_{\mathbf{i}'}\n$ which is positively equivalent to $\bbDelta_{\mathbf{i}}\n$.

We define the set of potentials on $G^*$, following notation of Section \ref{wedfodbc}
\[
  \Phi_{G^*}:=\left\{\frac{(E_i\cdot \Delta_{\omega_i,w_0\omega_i})_1}{(\Delta_{\omega_i,w_0\omega_i})_1}, \frac{(\Delta_{\omega_i,w_0\omega_i}\cdot E_{i'})_1}{(\Delta_{\omega_i,w_0\omega_i})_1}, \frac{((E_i\cdot \Delta_{\omega_i,w_0\omega_i})\circ\tau)_2}{(\Delta_{\omega_i,w_0\omega_i}\circ\tau)_2},\frac{( (\Delta_{\omega_i,w_0\omega_i}\cdot E_{i'})\circ \tau)_2}{(\Delta_{\omega_i,w_0\omega_i}\circ\tau)_2}\right\}_i,
\]
where as usual $i$ indexes the simple roots of $G$. Note that $\Phi_{G^*}\subset P(G^*,[\bbDelta_\mathbf{i}\n])$ by Proposition \ref{pro4.11}. Recalling Definition \ref{BKpotential}, we see that $\Phi_{G^*}$ comes from the Berenstein-Kazhdan potential $\Phi_{BK}$ on $G^{e,w_0}$. Note that $\Phi_{G^*}$ restricts to regular functions on the toric chart $\bbDelta_\mathbf{i}\n$, and so by Remark \ref{coneremark} the cone $\Ccal_{\Phi_{G^*}}(G^*,\bbDelta_\mathbf{i}\n)(\R)$ is convex.

\begin{theorem}\label{PPonG*}
  For any reduced word $\mathbf{i}$ of $w_0$, the quadruple $(G^*,\pi_{G^*},[\bbDelta_\mathbf{i}\n]_{\pi_{G^*}},\Phi_{G^*})$ is a positive Poisson variety.
\end{theorem}

\begin{proof}
  We must check that $\pi_{G^*}$ is weakly log-canonical in the coordinates given by $\bbDelta_\mathbf{i}$. Let $\Delta,\Delta'\in F(\mathbf{i})$. We must consider three types of brackets. Recall the notation of Definition \ref{projnotation}.
    
  (a) Bracket of type $\{\Delta_1,\Delta'_1\}_{G^*}$. By Proposition \ref{borelpoisson},
  \[
    \{\Delta_1,\Delta'_1\}_{G^*}=-\{\Delta,\Delta'\}_{B}\circ \pr_1.
  \]
  By Theorem 2.6 of \cite{KZ}, the bracket $\{\Delta\circ \xi,\Delta'\circ \xi\}_B$ is log-canonical on the open subset $G^{e,w_0}$, where $\xi:G^{e,w_0}\to G^{e,w_0}$ is the twist map of Proposition \ref{twistmap}. By Theorem 3.1 of \cite{GSV}, the twist map $\xi$ is anti-Poisson.
    
  (b) Bracket of type $\{(\Delta\circ\tau)_2,(\Delta'\circ\tau)_2\}_{G^*}$. Once we note $\tau:B_-\to B$ is anti-Poisson by Propostion \ref{tauAP}, the argument is the same as the previous case.
    
  (c) Bracket of type $\{\Delta_1,(\Delta'\circ\tau)_2\}_{G^*}$. By the definition of $\tau$ and (\ref{mixedbracket}), we have
  \begin{align}
     \{\Delta_1,(\Delta'\circ\tau)_2\}_{G^*} =\quad &   \label{simplifiedmixed} \frac{1}{2}\sum_{i=1}^{\operatorname{rank} \hhh} (X_i\cdot \Delta)_1( (X_i\cdot \Delta')\circ \tau)_2 - (\Delta\cdot X_i)_1((\Delta'\cdot X_i)\circ \tau)_2 \\
    +&\sum_{\alpha\in R^+} (E_\alpha\cdot \Delta)_1((E_\alpha\cdot \Delta')\circ \tau))_2 -(\Delta\cdot E_\alpha)_1((\Delta'\cdot E_\alpha)\circ \tau)_2.\nonumber
  \end{align}
  Write $\Delta=\Delta_{\omega_j,u\omega_j}$ and $\Delta'=\Delta_{\omega_k,v\omega_k}$. By (\ref{BonusIdentity}), the first sum
  \begin{align*}
    \frac{1}{2}\sum_{i=1}^{\operatorname{rank} \hhh} (X_i\cdot \Delta)_1( (X_i\cdot \Delta')\circ \tau)_2 - (\Delta\cdot X_i)_1((\Delta'\cdot X_i)\circ \tau)_2  \\
    =\frac{1}{2}\left( \sum_{i=1}^{\operatorname{rank} \hhh}  \omega_j(X_i) \omega_k(X_i)-u\omega_j(X_i) v\omega_k(X_i) \right)\Delta_1(\Delta'\circ \tau)_2  \\
    =\frac{1}{2} \left( K(\omega_j,\omega_k)-K(u\omega_j,v\omega_k)\right) \Delta_1 (\Delta'\circ \tau)_2 
  \end{align*}
  is log-canonical. For future reference, we note that here that the coefficient of $\Delta_1(\Delta'\circ \tau)_2$ is rational. We then write
  \begin{align}
    \{\Delta_1,(\Delta'\circ\tau)_2\}_{G^*} & = \Delta_1(\Delta'\circ\tau)_2\big( c + f),
  \end{align}
  where 
  \[
    f = \sum_{\alpha\in R^+} \frac{(E_\alpha\cdot \Delta)_1}{\Delta_1}\frac{((E_\alpha\cdot \Delta')\circ \tau)_2}{(\Delta'\circ \tau)_2} -\frac{(\Delta\cdot E_\alpha)_1}{\Delta_1} \frac{((\Delta'\cdot E_\alpha)\circ \tau)_2}{(\Delta'\circ \tau)_2}.
  \]
  By Theorem \ref{maintheorem} and Theorem \ref{maintheoremprime}, $f$ is weakly estimate dominated by $\Phi_{G^*}$, and thus weakly  dominated by $\Phi_{G^*}$.
\end{proof}

\begin{remark}
Theorem \ref{PPonG*} describes a number of toric charts (labeled by reduced words of $w_0$) making the dual Poisson algebraic group $G^*$ into a positive Poisson variety. In fact, there are many more charts with this property. In particular, every cluster chart on $G^{e,w_0}$ gives rise to a positive Poisson variety structure on $G^*$.  Note that the number of charts labeled by reduced words is always finite whereas the number of cluster charts is usually infinite.

It turns out that the proof of Theorem \ref{PPonG*} applies in the case of arbitrary cluster charts: Poisson brackets of the type $\{ f_1, f'_1\}_{G^*}$ reduce to Poisson brackets on $B$, and they are log-canonical in all cluster charts (by the standard properties of Poisson cluster varieties). And Poisson brackets of the type $\{ f_1, (f' \circ \tau)_2\}_{G^*}$ are weakly dominated by the potentials by Remark \ref{totally_positive}.

This viewpoint will be explored elsewhere. \end{remark}

\begin{example}
  \label{sl3pospoisson}
  Let $G={\rm SL}_3$ as in Examples \ref{sl3liegroupexample}, \ref{sl3factparexample}, and \ref{sl3posstrucexample}. The set of potentials on $SL_3^*$ is
\begin{align*}
\Phi_{G^*}= & \left\{ 
\frac{(\Delta_{13,23})_1}{(\Delta_{12,23})_1}, \frac{(\Delta_{2,3})_1}{(\Delta_{1,3})_1},\frac{(\Delta_{12,13})_1}{(\Delta_{12,23})_1}, \frac{(\Delta_{1,2})_1}{(\Delta_{1,3})_1},\right. \\ 
& \left. \frac{(\Delta_{13,23}\circ\tau )_2}{(\Delta_{12,23}\circ\tau)_2}, \frac{(\Delta_{2,3}\circ\tau)_2}{(\Delta_{1,3}\circ\tau)_2},\frac{(\Delta_{12,13}\circ\tau)_2}{(\Delta_{12,23}\circ\tau)_2}, \frac{(\Delta_{1,2}\circ\tau)_2}{(\Delta_{1,3}\circ\tau)_2}
\right\}.
\end{align*}
  From Theorem \ref{PPonG*}, we know that, with this choice of potentials $\Phi_{G^*}$, the chart $\bbDelta_\mathbf{i}\n$ is weakly log-canonical for $\pi_{{\rm SL}_3^*}$. As an example, we compute the bracket $\{(\Delta_{1,3})_1,(\Delta_{1,2}\circ \tau)_2\}$ in ${\rm SL}_3^*$; see (\ref{simplifiedmixed}) above:
\begin{align}
& \{(\Delta_{1,3})_1,(\Delta_{1,2}\circ \tau)_2\} \nonumber \\
&=\quad \frac{1}{2}\left(K(\omega_1,\omega_1)-K(w_0\omega_1,s_1\omega_1)\right) (\Delta_{1,3})_1(\Delta_{1,2} \circ \tau)_2 \nonumber  \\
  & \quad +\sum_{ij\in\{12,23,13\}} (E_{ij}\cdot \Delta_{1,3})_1((E_{ij}\cdot \Delta_{1,2})\circ \tau))_2 -(\Delta_{1,3}\cdot E_{ij})_1((\Delta_{1,2}\cdot E_{ij})\circ \tau)_2.\nonumber \\
  &  = \quad   c (\Delta_{1,3})_1( \Delta_{1,2}\circ \tau)_2   +(\Delta_{2,3})_1 (\Delta_{2,2}\circ\tau)_2+ (\Delta_{3,3})_1(\Delta_{3,2}\circ \tau)_2. \nonumber \\ 
  & = \quad (\Delta_{1,3})_1( \Delta_{1,2}\circ \tau)_2\left(
  c+\frac{(\Delta_{2,3})_1 (\Delta_{2,2}\circ\tau)_2}{(\Delta_{1,3})_1( \Delta_{1,2}\circ \tau)_2}+ \frac{(\Delta_{3,3})_1(\Delta_{3,2}\circ \tau)_2}{(\Delta_{1,3})_1( \Delta_{1,2}\circ \tau)_2}
  \right).
  \label{bracketexample}
\end{align}
One can check that the two terms $\frac{(\Delta_{2,3})_1 (\Delta_{2,2}\circ\tau)_2}{(\Delta_{1,3})_1( \Delta_{1,2}\circ \tau)_2}+ \frac{(\Delta_{3,3})_1(\Delta_{3,2}\circ \tau)_2}{(\Delta_{1,3})_1( \Delta_{1,2}\circ \tau)_2}$ of (\ref{bracketexample}) are indeed weakly estimate-dominated by $\Phi_{G^*}$.
\end{example}

\section{Tropicalization of Poisson Structures}

In this Section, we pass to complex points of our varieties, which have thus far been defined over $\Q$. Specifically, a \emph{complex positive variety} $X(\C)$ is the complex points of the positive variety $(X,[\theta])$. A toric chart $\theta:\G_m^n \to X$ induces an open embedding $(\C^\times)^n\to X(\C)$, which we also call a toric chart, and the definitions of the previous sections extend similarly under base change.

 In the remainder of this paper, we consider only complex positive  varieties and their real forms. For simplicity we write $X=X(\C)$. In particular, taking complex points of the Poisson algebraic groups of Section \ref{PLG} gives complex Poisson-Lie groups $G=G(\C)$ and $G^*=G^*(\C)$.
 
We define (under some extra assumptions) a tropicalization $\mathcal{C} \times \mathbb{T}$ of a complex positive Poisson variety $X$. Here $\mathcal{C}$ is a cone, $\mathbb{T}$ is a real torus, the tropicalization carries a constant Poisson structure. In particular, a dual Poisson-Lie group $G^*=G^*(\C)$ we obtain that $\mathcal{C}$ is the strict extended string cone.

\subsection{Real Forms of Poisson Structures}

\label{realforms1}

In this section we introduce real forms of holomorphic Poisson structures. Much of this material is well known; see for instance \cite{Xu}.

Let $(X,\pi)$ be a complex manifold with holomorphic Poisson structure $\pi\in \Lambda^2(T^{1,0} X)$. Let $\pi=\pi_R+i\pi_I$ be the decomposition of $\pi$ into real and imaginary parts; it is well-known that $\pi_R,\pi_I\in \Gamma(\Lambda^2(TX))$ are (real) Poisson bivectors. Let $\overline{\tau}:X\to X$ be an anti-holomorphic involution of $X$, which satisfies $\overline{\tau}(\pi_R)=\pi_R$. In this case, we say $\overline{\tau}$ is \emph{Poisson}. 

\begin{remark}  Equivalently, $\overline{\tau}$ is Poisson if $\overline{\tau}(\pi_I)=-\pi_I$. Extending $\overline{\tau}$ conjugate-linearly to $TX\otimes \C$, this is equivalent to the condition $\overline{\tau}(\pi)=\pi$. 
\end{remark} 

Let $K\subset \Fix(\overline{\tau})$ be a (real) open submanifold of the fixed points of $\overline{\tau}$. For any $p\in K$, decompose $T_pX$ as $T_pX=(T_pX)^{\overline{\tau}}+(T_pX)^{-\overline{\tau}}$, where
\[
	(T_pX)^{\overline{\tau}}=\{v\in T_pX:\overline{\tau}(v)=v\},\quad (T_pX)^{-\overline{\tau}}=\{v\in T_pX: \overline{\tau}(v)=-v\}. 
\]
As shown in \cite{Xu}, $\pi_R$ can be decomposed as $\pi_R(p)=\pi_R^{\overline{\tau}}(p)+\pi_R^{-\overline{\tau}}(p)$, where
\[
	\pi_R^{\overline{\tau}}(p)\in \Lambda^2 (T_pX)^{\overline{\tau}},\quad\pi_I^{-\overline{\tau}}(p)\in \Lambda^2 (T_pX)^{-\overline{\tau}}.
\]
\begin{lemma}[\cite{Xu}]
	Using the notation above, $\pi_R^{\overline{\tau}}$ is a Poisson bivector on $K$.
\end{lemma}


The pair $(K, \pi_R^{\overline{\tau}})$ is called a \emph{real form} of $X$. 

\begin{example}
  Consider any holomorphic Poisson structure
  \[
    \pi=\sum_{i,j} \pi_{ij}(z)\partial_{z_i}\wedge\partial_{z_j}
  \]
  on $\C^n$, where $\pi_{ij}(z)$ are holomorphic functions. Let $\overline{\tau}$ be the anti-holomorphic involution of $\C^n$ given by $\overline{\tau}(z)=\overline{z}$. Then the set of fixed points of $\overline{\tau}$ is $K=\R^n\subset \C^n$. Thus $\overline{\tau}(\pi_R)=\pi_R$ if and only if $\pi_{ij}(\overline{z})=\overline{\pi_{ij}(z)}$. Write
  \[
    \partial_{z_i}=\frac{1}{2}(\partial_{x_i}-\sqrt{-1}\partial_{y_i}),
  \]
  and direct calculation shows
  \begin{equation}\label{realocal}
    \pi_R^{\overline{\tau}}=\frac{1}{4}\sum_{i,j} \pi_{ij}(x)\partial_{x_i}\wedge\partial_{x_j}.
  \end{equation}
\end{example}

\begin{lemma}
\label{realocallemma}
  Let $X$ be a complex manifold with holomorphic Poisson structure $\pi$. Let $\overline{\tau}$ be an anti-holomorphic Poisson involution of $(X,\pi)$. Let $(K,\pi_{R}^{\overline{\tau}})$ be the corresponding real form. Then we have
  \[
    \{f_1|_K,f_2|_K\}_{\pi_R^{\overline{\tau}}}=\frac{1}{4}\{f_1,f_2\}_{\pi}|_{K},
  \]
  where the $f_i$ are  holomorphic functions on an open subset $U\subset X$ satisfying $f_i(\overline{\tau}(z))=\overline{f_i(z)}$.
\end{lemma}

\begin{proof}
  We only need to show the lemma in a neighborhood of each fixed point of $\overline{\tau}$. Choose holomorphic local coordinates $z_1,\dots,z_n$ such that $\overline{\tau}$ is given by $\overline{\tau}(z)=\overline{z}$. Then in these coordinates, the $f_i$ satisfy $f_i(\overline{z})=\overline{f_i(z)}$. Set $z_j=x_j+\sqrt{-1}y_j$ and let $x=(x_1,\dots,x_n)$ and and $y=(y_1,\dots,y_n)$. Write $f_i(z)=u_i(x,y)+\sqrt{-1}v_i(x,y)$, where $u_i,v_i\in C^\infty(X)$ are smooth real-valued functions on $X$. Since $f_i(\overline{z})=\overline{f_i(z)}$, we know $u_i(x,y)=u_i(x,-y)$ and $v_i(x,y)=-v_i(x,-y)$. Thus $\partial_{y_j}u_i|_{y=0}=0$. Then by the Cauchy–Riemann equations, we have
  \begin{align*}
    \Big(\partial_{z_j}f_i(z)\Big)\Big|_{K}&=\frac{1}{2}\Big(\partial_{x_j}u_i+\partial_{y_j}v_i-\sqrt{-1}\partial_{y_j}u_i+\sqrt{-1}\partial_{x_j}v_i\Big)\Big|_{y=0}\\
    &=\left. \frac{1}{2}\Big(\partial_{x_j}u_i+\partial_{y_j}v_i\Big)\right|_{y=0}=\partial_{x_j}f_i(x)
  \end{align*}
  By Equation \eqref{realocal}, we get the conclusion.
\end{proof}

\begin{prop} 
\label{realbracketisok} 
  Let $X=(\C^\times)^n$ with holomorphic Poisson bivector $\pi$, let $\tau:X\to X$ be an algebraic involution of $X$, and let $\overline{\tau}:X\to X$ be the anti-holomorphic involution given by $\overline{\tau}(z):=\tau(\bar{z})=\overline{\tau(z)}$. Assume $\overline{\tau}$ is Poisson, and consider the real form $(K,\pi_{R}^{\overline{\tau}})$. Let $f_1,f_2$ be holomorphic functions on an open subset $U\subset X$. Then
  \[
    \{f_1|_K,f_2|_K\}_{\pi_R^{\overline{\tau}}}=\frac{1}{4}\{f_1,f_2\}_{\pi}|_{K}.
  \]
\end{prop} 

\begin{proof}
  Let $g_i:= f_i+f_i\circ \tau$ and $h_i:=f_i-f_i\circ \tau$. Then $g_i$ and $\sqrt{-1}h_i$ satisfy the condition from Lemma \ref{realocallemma}: 
  \[ g_i(\overline{\tau}(z))=\overline{g_i(z)},\qquad \sqrt{-1}h_i(\overline{\tau}(z))=\overline{\sqrt{-1}h_i(z)}.\] We then compute:
  \begin{align*}
  \{f_1|_K,f_2|_K\}_{\pi_R^{\overline{\tau}}} & = \frac{1}{4}\{g_1|_K+h_1|_K, g_2|_K+h_2|_K\}_{\pi_R^{\overline{\tau}}} \\ 
  & = \frac{1}{4}\Big(
  \{g_1|_K,g_2|_K\}_{\pi_R^{\overline{\tau}}}-\sqrt{-1}\{g_1|_K,\sqrt{-1}h_2|_K\}_{\pi_R^{\overline{\tau}}} \\
  & \quad -\sqrt{-1}\{\sqrt{-1}h_1|_K,g_2|_K\}_{\pi_R^{\overline{\tau}}}-\{\sqrt{-1}h_1|_K,\sqrt{-1}h_2|_K\}_{\pi_R^{\overline{\tau}}}\Big)
  \\ 
  & = \frac{1}{16}\Big(
  \{g_1,g_2\}_\pi|_K-\sqrt{-1}\{g_1,\sqrt{-1}h_2\}_\pi|_K \\
  & \quad -\sqrt{-1}\{\sqrt{-1}h_1,g_2\}_\pi|_K-\{\sqrt{-1}h_1,\sqrt{-1}h_2\}_\pi|_K
  \Big) \\ 
  & = \frac{1}{16}\{g_1+h_1,g_2+h_2\}_\pi|_K = \frac{1}{4}\{f_1,f_2\}_\pi|_K.
  \end{align*}
\end{proof}

\subsection{Real Forms of Positive Poisson Varieties}

In this section we define real forms of positive Poisson varieties, and give as an example the real form $K^*$ of the positive Poisson variety $G^*$.

Let $(X,[\theta])$ be a complex positive variety, with $\theta:(\C^\times)^n\to X$. Then complex conjugation $\overline{(\cdot)}:(\C^\times)^n\to (\C^\times)^n$ defines an anti-holomorphic involution on the open subvariety $\theta((\C^\times)^n)\subset X$. Since the transition maps between charts in $[\theta]$ are positive, their components have rational coefficients and they commute with complex conjugation of $(\C^\times)^n$. So $\overline{(\cdot)}$ extends to the open subvariety $\bigcup_{\theta'\in[\theta]}\theta'((\C^\times)^n)$. In particular, if the charts in $[\theta]$ cover $X$, the involution $\overline{(\cdot)}$ is defined on all of $X$.

\begin{definition}\label{posrealformdef}
  Let $(X,\pi,[\theta]_{\pi},\Phi)$ be a complex positive Poisson variety, with $\theta:(\C^\times)^n\to X$ a $\pi$-compatible chart. Let $\tau:X\to X$ be an involution, which is also a map of positive varieties with potential. Then $\overline{\tau}:X\to X$ is an anti-holomorphic involution of $X$, where $\overline{\tau}(x):=\tau(\bar{x})=\overline{\tau(x)}$. Assume $\overline{\tau}$ is Poisson.

  
  We impose the additional conditions:
  \begin{enumerate}[label=(\alph*)]
    \item In the standard coordinates $z_1,\dots,z_n$ on $(\C^\times)^n$, the positive transformation $\tau$ is monomial.
    \item The log-canonical part $\pi_{ij}$ of the bracket of coordinate functions $\{z_i,z_j\}$ is pure imaginary.
  \end{enumerate}
  In this case we call the tuple $(X,\pi,\theta,\Phi,\tau)$ a \emph{framed positive Poisson variety with real form}. 
  Maps from $(X_1,\pi_1,\Theta_{\pi_1},\Phi_1,\tau_1)$ to $(X_2,\pi_2,\Theta_{\pi_2},\Phi_2,\tau_2)$ are maps of complex positive Poisson varieties which intertwine $\tau_1$ and $\tau_2$. The corresponding category is denoted $\mathbf{PosPoiss}^\bullet_\R$.
\tria \end{definition}

\begin{theorem} 
\label{realformofG^*}
  Let $(G^*,\pi_{G^*})$ be the dual Poisson-Lie group given by complex points of the Poisson algebraic group $G^*$ introduced in Section \ref{PLG}, and let $\tau$ the involution as in Proposition \ref{tauAP}:
  \[
    \tau:G^*\to G^*\ :\  (b^+,b^-)\mapsto (((b^-)^T)\n,((b^+)^T)\n).
  \]
  
  Then for a reduced word $\mathbf{i}$ for $w_0$, the tuple $(G^*,\sqrt{-1}\pi_{G^*},\bbDelta_\mathbf{i}\n,\Phi_{G^*}, \tau)$ is a framed positive Poisson variety with real form. 
\end{theorem}
  
\begin{proof}
  From the expression of the $r$-matrix $r_D$, we see that $\tau$ is a holomorphic anti-Poisson involution of $(G^*,\pi_{G^*})$. It follows that  $\overline{\tau}$ is an anti-holomorphic Poisson involution of $(G^*,\sqrt{-1}\pi_{G^*})$. In the toric coordinates given by $\bbDelta_\mathbf{i}\n$, the involution $\tau$ simply permutes the coordinates and thus $\tau$ satisfies condition (a) of Definition \ref{posrealformdef}. Also, $\tau^*\Phi_{G^*}=\Phi_{G^*}$, so $\tau$ is a map of positive varieties with potential. By the proof of Theorem \ref{PPonG*}, the weakly log-canonical part of the bracket $\sqrt{-1}\pi_{G^*}$ is a rational multiple of $\sqrt{-1}$, and therefore satisfies condition (b) of Definition \ref{posrealformdef}.
\end{proof}

We next want to specify a canonical choice for the real form of $(X,\pi,\theta,\Phi,\tau)$; note the real form $K$ was not uniquely determined in Section \ref{realforms1}. We first introduce some notation.

\begin{definition}[Notation]
\label{notation1forPT}
Let $(X,\pi,\theta:(\C^\times)^n\to X,\Phi,\tau)$ be a framed positive Poisson variety with real form. We set $(X,\theta)^t(\R):=(X,\theta)^t\otimes_\Z \R$ to be the extension by scalars. For simplicity, sometimes we use $X^t(\R)$ for $(X,\theta)^t(\R)$. Let
\[
L_{X,\theta,\tau}:=\{\xi\in X^t(\R)| \tau^t(\xi)=\xi\}
\]
be the fixed points of the tropical involution $\tau^t$, and let
\[
\T_{X,\theta,\tau}:=\{g\in (S^1)^n\subset (\C^\times)^n| \tau(\bar{g})=g\}_0
\]
be the identity component of the fixed points of the anti-holomorphic involution $\overline{\tau}|_{(S_1)^n}$. Sometimes we write $L_X=L_{X,\theta,\tau}$ and $\T_X=\T_{X,\theta,\tau}$ for brevity. Finally, let $s>0$ be a real parameter. Define the map
\begin{align*}
    E_{X,\theta,s}: (X,\theta)^t(\R)\times (S^1)^n & \to (\C^\times)^n; \\
    (\xi_1,\dots,\xi_n,e^{\sqrt{-1}\nu_1},\dots,e^{\sqrt{-1}\nu_n}) & \mapsto (e^{s\xi_1+\sqrt{-1} \nu_1},\dots, e^{s\xi_n+\sqrt{-1}\nu_n}).
\end{align*}
Recall our assumption that $\tau$ is a monomial transformation in the toric coordinates given by $\theta$. Thus $\overline{\tau}$ preserves the compact torus $(S^1)^n\subset (\C^\times)^n$, and $L_{X,\theta,\tau} \subset X^t(\R)$ is a linear subspace. \tria \end{definition}

Pulling back $\overline{\tau}$ by $\theta\circ E_{X,\theta,s}$, we find:

\begin{prop}
\label{realformasimage}
Let $(X,\pi,\theta,\Phi,\tau)$ be a framed positive Poisson variety with real form. Then for all $s,s'>0$, we have 
\[
(\theta\circ E_{X,\theta,s})(L_X\times\T_X)=(\theta\circ E_{X,\theta,s'})(L_X\times\T_X)\subset \Fix(\overline{\tau})\cap \theta((\C^\times)^n) \subset X.
\]
\end{prop}

\begin{definition}
    \label{Krealformdefn} Let $(X,\pi,\theta,\Phi,\tau)$ be a framed positive Poisson variety with real form. Then the \emph{real form} of $(X,\pi,\theta,\Phi,\tau)$ is $(\Re(X,\pi,\theta,\Phi,\tau),\pi_{\Re(X,\theta,\Phi,\tau)})$, where 
    \[
    \Re(X,\pi,\theta,\Phi,\tau)=(\theta\circ E_{X,\theta,s'})(L_X\times\T_X)\subset  \Fix(\overline{\tau})\cap \theta((\C^\times)^n),\quad \pi_{\Re(X,\theta,\Phi,\tau)}=4\pi_R^{\overline{\tau}}.
    \]
    We usually write $(\Re(X),\pi_{\Re(X)})$ for simplicity.
    Note that the definition does not depend on the choice of $s>0$, and that, (up to a scalar multiple of the bracket), the Poisson manifold $(\Re(X),\pi_{\Re(X)})$ is a real form in the sense of the previous section. In summary, we have the following diagram
    \begin{equation*}
    \begin{tikzcd}[column sep=4.3em]
        (X,\theta)^t(\R)\times (S^1)^n \arrow{r}{E_{X,\theta,s}}[swap]{\sim} &(\C^\times)^n \arrow[hookrightarrow]{r}{\theta} & X \\
        L_{X,\theta,\tau}\times \T_{X,\theta,\tau} \arrow{rr}{\theta\circ E_{X,\theta,s}}[swap]{\sim} \arrow[hookrightarrow]{u} & & \Re(X). \arrow[hookrightarrow]{u}
    \end{tikzcd}
    \end{equation*}
\tria \end{definition}

\begin{remark}
\label{posrealformremark}
 By Proposition \ref{realbracketisok}, we have $\{z_i|_{\Re(X)},z_j|_{\Re(X)}\}_{\pi_R^{\overline{\tau}}}=\frac{1}{4}\{z_i,z_j\}_\pi|_{\Re(X)}$. 
The factor of $4$ in the definition of $\pi_{\Re(X)}$ is simply to make formulas simpler.
\end{remark}

\begin{prop}
\label{mapsofrealforms}
Let $
f:(X_1,\pi_1,\theta_1,\Phi_1,\tau_1)\to (X_2,\pi_2,\theta_2,\Phi_2,\tau_2)
$ be a map of framed positive Poisson varieties with real forms. Then $\theta_2\n\circ f\circ\theta_1$ restricts to a Poisson map \[ (\theta_1\n(\Re(X_1)),\theta_1^*\pi_{\Re(X_1)})\to (\theta_2\n(\Re(X_2)),\theta_2^*\pi_{\Re(X_2)})\] whenever it is defined.
\end{prop}

\begin{proof}
Since $f$ is a map of framed positive Poisson varieties with real forms, we have $f\circ \overline{\tau_1}=\overline{\tau_2}\circ f$. Because $f$ is a positive rational map, it takes $\theta_1(\T_{X_1})$ to $\theta_2(\T_{X_2})$, and so it takes the component $\Re(X_1)$ of $\Fix(\overline{\tau}_1)$ to the component $\Re(X_2)$ of $\Fix(\overline{\tau}_2)$. That the map is Poisson follows from Proposition \ref{realbracketisok}.
\end{proof}

\begin{remark} 
\label{realformG^*remark} Let $K^*\subset G^*$ be the identity component of the fixed points of $\overline{\tau}$. From Proposition \ref{realformofG^*} we know $(G^*,\sqrt{-1}\pi_{G^*},\bbDelta_\mathbf{i}\n,\Phi_{G^*}, \tau)$ is a positive Poisson variety with real form. The real form given by Definition \ref{Krealformdefn} is an open dense submanifold of $K^*$, and its Poisson structure extends to the Poisson structure $\pi_{K^*}=4(\sqrt{-1} \pi_{G^*})_R^{\overline{\tau}}$ on $K^*$. We then think of the real form of $(G^*,\sqrt{-1}\pi_{G^*},\bbDelta_\mathbf{i}\n,\Phi_{G^*}, \tau)$ as a coordinate neighborhood on $K^*$.

The Poisson-Lie group $(G,\pi_G)$ has a compact real form $(K,\pi_K)$; see \cite{AD}. It is shown, see for example \cite{ES,Xu}, that $(K^*,\pi_{K^*})$ is a dual Poisson-Lie group of $(K,\pi_K)$. From the Iwasawa decomposition $G=KAN$ of $G$, we may identify $K^*=AN$.
 \end{remark}
 
 \begin{example}
   \label{sl3realformexample}
   Let $G={\rm SL}_3(\C)$ as in Examples \ref{sl3liegroupexample}, \ref{sl3factparexample}, \ref{sl3posstrucexample}, and \ref{sl3pospoisson} above. It is easy to see that $\tau$ is a positive involution of $({\rm SL}_3^*,[\bbDelta_\mathbf{i}\n])$, and $\tau^t$ preserves $\Ccal_{\Phi_{G^*}}$. In the toric coordinates $z_1,\dots,z_8$ given by $\bbDelta_\mathbf{i}\n$, the anti-holomorphic involution $\overline{\tau}: (\C^\times)^8\to (\C^\times)^8$ becomes,
\[ 
\overline{\tau}(z_1,\dots,z_8)=(\bar{z}_4,\bar{z}_5,\bar{z}_6,\bar{z}_1,\bar{z}_2,\bar{z}_3,\bar{z}_7,\bar{z}_8).
\]
Fixed points are those of the form 
\[ 
(z_1,z_2,z_3,\bar{z}_1,\bar{z}_2,\bar{z}_3,x_7,x_8),\qquad\text{ where }x_7,x_8\in \R\setminus\{0\}.
\]
Projecting parallel to the fourth, fifth, and sixth coordinates and restricting to the positive part of the last two coordinates gives coordinates on the usual presentation of $K^*$ as the group of complex upper-triangular matrices, with positive real entries along the diagonal, as in \cite{Anton, AM}.
 \end{example}

\subsection{Partial Tropicalization} 

In this section we construct for each framed positive Poisson variety with real form $(X,\pi,\theta,\Phi,\tau)$ a real manifold $PT(X,\pi,\theta,\Phi,\tau)$ with constant Poisson bracket $\pi_{PT}$ called its \emph{partial tropicalization}. We extend partial tropicalization to a functor
\[
PT:\mathbf{PosPoiss^\bullet_\R}\to \mathbf{PTrop}
\]
to the category of partial tropicalizations, defined below. We show in Theorem \ref{scalingtrans} below that partial tropicalization can be thought of as the limit of a 1-parameter family of coordinates on the real form $\Re(X)$.

\begin{definition}[Notation]
\label{realnotation}
Let $(X,\theta_X:(\C^\times)^n\to X,\Phi_X)$ 
be a framed complex positive variety with potential. 
Let \[ \mathcal{C}_\Phi(X,\theta_X)(\mathbb{R}) = \{\xi\in (X,\theta_X)^t(\R)|\Phi^t(\xi)<0\};
\]
see also Remark \ref{coneremark}. For simplicity, sometimes we use 
$\Ccal_{\Phi}(\R)$ for 
$\mathcal{C}_\Phi(X,\theta_X)(\mathbb{R}).$ 
\tria \end{definition}

\begin{definition}
\label{partropdef}
Let $(X,\pi,\theta,\Phi,\tau)$ be a framed positive Poisson variety with real form. Define
\[
PT(X,\pi,\theta,\Phi,\tau):=(\Ccal_\Phi(\R)\cap L_{X,\theta,\tau})\times \T_{X,\theta,\tau},
\]
where $L_{X,\theta,\tau}$ and $\T_{X,\theta,\tau}$ are given in Definition \ref{notation1forPT}.


Let $z_1,\dots,z_n$ be the standard coordinate functions on the toric chart $\theta:(\C^\times)^n\to X$. Consider induced involution $\tau^*$ of the field of rational functions $\C(z_1,\dots,z_n)$. Recall from Proposition \ref{monomiallogcan} that $\{z_i,\tau^*(z_j)\}$ is weakly log-canonical. Denote by $\pi_{i\tau(j)}$ the log-canonical part of the bracket, and note that by assumption (b) of Definition \ref{posrealformdef} and the proof of Proposition \ref{monomiallogcan}, we have $\pi_{i,\tau(j)}\in \sqrt{-1}\R$ is pure imaginary. Also, recall that the coordinates $z_1,\dots,z_n$ on $(\C^\times)^n$ induce coordinates $\xi_1,\dots,\xi_n$ on $X^t(\R)$.
Writing $(e^{\sqrt{-1}\nu_1},\dots,e^{\sqrt{-1}\nu_n})\in (S^1)^n$ gives us local coordinates $\nu_1,\dots, \nu_n$ on $(S^1)^n\subset(\C^\times)^n$. We define a constant real bracket on $X^t(\R)\times (S^1)^n$ as follows:
\begin{align*}
    \{\xi_i,\nu_j\} &  := \frac{-\sqrt{-1}}{2}(\pi_{ij}-\pi_{i\tau(j)}), \\
    \{\nu_i,\xi_j\} & := \frac{-\sqrt{-1}}{2}(\pi_{ij}+\pi_{i\tau(j)}), \\
    \{\xi_i,\xi_j\} & := 0, \\
    \{\nu_i,\nu_j\} & := 0.
\end{align*}
We set $\pi_{PT}$ to be the restriction of this bracket to $PT(X,\pi,\theta,\Phi,\tau)\subset X^t(\R)\times (S^1)^n$. The pair $(PT(X,\pi,\theta,\Phi,\tau),\pi_{PT})$ is the \emph{partial tropicalization} of $(X,\pi,\theta,\Phi,\tau)$.
\tria \end{definition}

 We will see $\pi_{PT}$ is skew-symmetric and well defined after giving a scaling interpretation for partial tropicalization in Theorem \ref{scalingtrans}.


\begin{lemma}
  \label{weakdomscaling} Let $f\in \tilde{P}_\Phi(X,[\theta])$ be a function which is weakly dominated by $\Phi$, and assume $f$ restricts to a regular function on the toric chart $\theta:(\C^\times)^n\to X$. Then, for every $\delta>0$ and every point
  in $\Ccal_{\Phi}(\R)\times (S^1)^{n}$ the function $f\circ\theta\circ  E_{X,\theta,s}$
  is bounded
  \[
    |f\circ\theta\circ E_{X,\theta,s}|\le e^{-s\delta}
  \]
  for $s$ sufficiently large.
\end{lemma}

\begin{proof}
    It suffices to show the estimate for $f\in P_\Phi(X,[\theta])$. By restricting to our toric chart we may assume $(X,\theta)=((\C^\times)^n,\Id)$. Then by assumption $f\in \C[z_1^{\pm1},\dots, z_n^{\pm1}]$ is regular on $(\C^\times)^n$, and so by the triangle inequality we may assume $\nu_i=0$ for all $i$. We write $(f\circ E_s)(\xi_1,\dots,\xi_n)$ for $(f\circ E_{(\C^\times)^n,\Id,s})(\xi_1,\dots,\xi_n,1,\dots,1)$. 
    Now, by a standard bounding argument as in \cite{Mi}, we get 
    \[
      \lim_{s\to \infty} \frac{1}{s} \log (f\circ E_s)(\xi_1,\dots,\xi_n)=f^t(\xi_1,\dots,\xi_n).
    \]
    For $(\xi_1,\dots,\xi_n)\in \Ccal_\Phi$, we have $f^t(\xi_1,\dots,\xi_n)<0$ since $f$ is dominated. Then for sufficiently large $s$, 
    \[
      \frac{1}{s} \log (f\circ E_s)(\xi_1,\dots,\xi_n)\le -\delta<0.
    \]
    Therefore, since exponential is monotonic,
    \[
      0<e^{s\frac{1}{s} \log (f\circ E_s)(\xi_1,\dots,\xi_n)} \le e^{-s\delta}.
    \]
    Thus, $(f\circ E_s)(\xi_1,\dots,\xi_n)=e^{s\frac{1}{s}\log (f\circ E_s)(\xi_1,\dots,\xi_n)}\to 0$ exponentially quickly as $s\to \infty$.
\end{proof}

 We introduce the scaled Poisson bivector $\pi_s:=s\pi_{\Re(X)}$ on the real form $\Re(X)=(\theta\circ E_{X,\theta,s})(L_{X,\theta,\tau}\times\T_{X,\theta,\tau})$.

\begin{theorem}[Partial Tropicalization as a Limit]\label{scalingtrans}
  Under the change of coordinates $\theta\circ E_{X,\theta,s}$, in the limit $s\to \infty$, the bivector $\pi_s$ converges to the constant bivector $\pi_{PT}$ on $PT(X,\pi,\theta,\Phi,\tau)\subset L_{X,\theta,\tau}\times\T_{X,\theta,\tau}$.
\end{theorem}

\begin{proof}
     By restricting to our toric chart, we may assume $(X,\theta)=((\C^\times)^n,\Id)$. Note that since $\theta$ is an open embedding, regular functions on $X$ restrict to Laurent polynomials in the coordinates $z_i$ on $(\C^\times)^n$.
     We compute:
    \begin{align}
      \{z_i,z_j\}_s & = \{e^{s\xi_i+\sqrt{-1}\nu_i},e^{s\xi_j+\sqrt{-1}\nu_j}\}_s \nonumber \\
      & =se^{s\xi_i+\sqrt{-1}\nu_i} e^{s\xi_j+\sqrt{-1}\nu_j}(\pi_{ij}+f_{ij}), \label{brack1}
    \end{align}
    where we write $f_{ij}=f_{ij}(e^{s\xi_1+\sqrt{-1} \nu_1},\dots, e^{s\xi_n+\sqrt{-1}\nu_n})$.
    On the other hand, $\{\cdot,\cdot \}_s$ is a biderivation, so
    \begin{align}\label{brack2}
      \{e^{s\xi_i+\sqrt{-1}\nu_i},e^{s\xi_j+\sqrt{-1}\nu_j}\}_s =\ \  &e^{s\xi_i+\sqrt{-1}\nu_i} e^{s\xi_j+\sqrt{-1}\nu_j}\Big(s^2\{\xi_i,\xi_j\}_s \\ 
      +&\sqrt{-1}s(\{\xi_i,\nu_j\}_s+\{\nu_i,\xi_j\}_s)-\{\nu_i,\nu_j\}\Big). \nonumber
    \end{align}
    Combining (\ref{brack1}) and (\ref{brack2}) gives
    \begin{equation}\label{brack3}
      s^2\{\xi_i,\xi_j\}_s+\sqrt{-1}s(\{\xi_i,\nu_j\}_s+\{\nu_i,\xi_j\}_s)-\{\nu_i,\nu_j\}_s= s(\pi_{ij}+f_{ij}).
    \end{equation}
    The condition that $\pi_s$ is a real bivector on $\Re(X)$ implies that for complex-valued functions $f,g\in C^\infty(\Re(X),\C)$, the bracket respects complex conjugation  $\{\overline{f},\overline{g}\}_s=\overline{\{f,g\}_s}$. So, similarly to (\ref{brack3}) we find from considering the bracket $\{e^{s\xi_i-\sqrt{-1}\nu_i},e^{s\xi_j-\sqrt{-1}\nu_j}\}_s$, that
    \begin{equation}\label{brack4}
      s^2\{\xi_i,\xi_j\}_s-\sqrt{-1}s(\{\xi_i,\nu_j\}_s+\{\nu_i,\xi_j\}_s)-\{\nu_i,\nu_j\}_s= s(\overline{\pi}_{ij}+\overline{f}_{ij}).
    \end{equation}
    Recall we have assumed the log-canonical part $\pi_{ij}$ of the bracket $\{z_i,z_j\}$ is pure imaginary.
    Putting together (\ref{brack3}) and (\ref{brack4}), we get
    \begin{align} \label{brack5}
      s^2\{\xi_i,\xi_j\}_s-\{\nu_i,\nu_j\}_s & =s(\text{w.d.t.}), \\ \label{brack6}
      \{\xi_i,\nu_j\}_s+\{\nu_i,\xi_j\}_s & = -\sqrt{-1}\pi_{ij}+\text{w.d.t.},
    \end{align}
    where $\text{w.d.t.}$ stands for weakly dominated terms.
    
    Restricting to the fixed locus of $\overline{\tau}$, we have the relation $\overline{z}_i=\tau^* (z_j)$. Therefore, on $\Re(X) \subset \Fix(\overline{\tau})$ we have
    \begin{equation}\label{brack7}
      \{z_i,\overline{z}_j\}_{\Re(X)}=\{z_i,\tau^*(z_j)\}_{\Re(X)}=z_i\tau^*(z_j)(\pi_{i\tau(j)}+f_{i\tau(j)})=z_i\overline{z}_j(\pi_{i\tau(j)}+f_{i\tau(j)}).
    \end{equation}
    Repeating calculations similar to those before (\ref{brack5}) and (\ref{brack6}) gives
    \begin{align}\label{brack8}
      s^2\{\xi_i,\xi_j\}_s+\{\nu_i,\nu_j\}_s & =s(\text{w.d.t.}), \\\label{brack9}
      -\{\xi_i,\nu_j\}_s+\{\nu_i,\xi_j\}_s & = -\sqrt{-1}\pi_{i\tau(j)} +\text{w.d.t.}
    \end{align}
     Combining (\ref{brack5}), (\ref{brack6}), (\ref{brack8}), and (\ref{brack9}), and applying Lemma \ref{weakdomscaling} gives the result.
\end{proof}

\begin{cor}\label{tropbracketisok}
  The constant bracket $\pi_{PT}$ is well defined on $PT(X,\pi,\theta,\Phi,\tau)$. It is skew-symmetric, and hence Poisson.
\end{cor}

\begin{proof} 
  Follows from Theorem \ref{scalingtrans} and the fact that  $\Lambda^2(T_p \Re(X))$ is a (closed) linear subspace of $T_p \Re(X) \otimes T_p \Re(X)$, for any $p\in \Re(X)$.
\end{proof}

\begin{definition}\label{partialtrop}
  We define the category $\mathbf{PTrop}$ as follows. Objects are pairs $(M,\pi)$, where $M=\Ccal\times \T\subset \R^n\times (S^1)^n$ is the product of an open real cone and a subtorus of the compact torus, and $\pi$ is a constant Poisson bivector on $M$ which restricts to $0$ on $\Ccal$ and $\T$. We consider $\R^n=\Z^n\otimes_\Z \R$ to be the extension of the integer lattice, and have the exponential map $\exp(2\pi\sqrt{-1}\cdot):\R^n\to (S^1)^n$.
  Arrows in $\mathbf{PTrop}$ from $\Ccal\times \T$ to $ \Ccal'\times\T'\subset \R^m\times (S^1)^m$ are continuous piecewise linear maps $f:\R^n\to \R^m$ which are homogeneous in the sense of (\ref{homogeq}), have $f(\Z^n)\subset \Z^m$, and are subject to the following conditions.
  \begin{enumerate}
    \item $f(\Ccal)\subset \Ccal'$ .
    \item On each open linearity chamber $C\subset \R^n$ of $f$, note that there is a unique Lie group homomorphism $\exp(f|_C):(S^1)^n\to (S^1)^m$ induced by $\exp(2\pi \sqrt{-1}\cdot)$.
    When $C\cap \Ccal\ne \emptyset$, we require that the induced maps $\exp(f|_C):(S^1)^n\to (S^1)^m$ restrict to maps of the subgroups $\T\to \T'$.
    \item On each open linearity chamber $C$, the restriction of the map $f\times \exp(f|_C):(C\cap \Ccal)\times \T\to \Ccal'\times \T'$ is Poisson. 
  \end{enumerate}
  Note that an arrow in $\mathbf{PTrop}$ induces a Poisson map $(\Ccal\times \T, \pi)\to (\Ccal'\times\T',\pi')$, defined on an open subset of the domain.
\tria \end{definition}

\begin{definition}
  Let $f:(X,\pi_X,\theta_X,\Phi_X,\tau_X)\to(Y,\pi_Y,\theta_Y,\Phi_Y,\tau_Y)$ be a map of framed positive Poisson varieties with real form. Define the partial tropicalization $PT(f)$ of $f$ to be the continuous piecewise linear map $f^t:X^t(\R)\to Y^t(\R)$.
\tria \end{definition}

After proving the following Lemma, we will show in Theorem \ref{parttropfunctor} that the assignment $PT$ is functorial.

\begin{lemma}\label{scalingandtrop}
  Let $(X,\theta:(\C^\times)^n\to X)$ be a framed complex positive variety. Let $f:(\C^\times)^n\to \C^\times$ be a positive rational function, and let $C\subset (X,\theta)^t(\R)$ be an open linearity chamber of $f^t$. Then for $f_s=E_{\C^\times,\Id,s}\n\circ f\circ E_{X,\theta,s}$, we have
  \begin{equation}\label{scalingandtropequation}
    \lim_{s\to \infty} f_s|_C= f^t|_C\times \exp(f^t|_C):C\times(S^1)^n\to \R\times S^1.
  \end{equation}
  In particular, for a point $(\xi,e^{\sqrt{-1}\nu})\in C\times (S^1)^n$, for large enough $s$ the value of $f_s(\xi,e^{\sqrt{-1}\nu})$ is always defined.
\end{lemma}

\begin{proof}
  For  $\xi:=(\xi_1,\dots,\xi_n)\in C$ and $e^{\sqrt{-1}\nu}:=(e^{\sqrt{-1}\nu_1},\dots,e^{\sqrt{-1}\nu_n})\in (S^1)^n$ we have
  \[
    f_s(\xi,e^{\sqrt{-1}\nu})=\left(\frac{1}{s}\log |f(e^{s\xi+\sqrt{-1}\nu})|, \frac{f(e^{s\xi+\sqrt{-1}\nu})}{|f(e^{s\xi+\sqrt{-1}\nu})|}\right).
  \]
  Write $f=A/B$, where $A$ and $B$ are positive polynomials in $z_1,\dots,z_n$. Then
  \begin{align*}
      A(e^{s\xi+\sqrt{-1}\nu})  =\sum_j e^{sM_j(\xi)+ \sqrt{-1} M_j(\nu)},\quad
      B(e^{s\xi+\sqrt{-1}\nu})  =\sum_k e^{sN_k(\xi)+ \sqrt{-1} N_k(\nu)},
  \end{align*}
  where $M_j$ and $N_k$ are linear polynomials with positive integer coefficients. Without loss of generality, assume $M_1(\xi)>M_j(\xi)$ for $j>1$ and $N_1(\xi)>N_k(\xi)$ for $k>1$. Note this uses that $\xi$ is in an open linearity chamber of $f^t$. Also, note that $f^t|_C$ is given by the linear polynomial $M_1-N_1$.

  We then have
  \begin{align*}
      \lim_{s\to\infty} \frac{A(e^{s\xi+\sqrt{-1}\nu})}{e^{sM_1(\xi)+ \sqrt{-1} M_1(\nu)}} =1, \quad
      \lim_{s\to\infty} \frac{B(e^{s\xi+\sqrt{-1}\nu})}{e^{sN_1(\xi)+ \sqrt{-1} N_1(\nu)}} =1,
  \end{align*}
  and so, dividing we get
  \begin{equation}\label{fraclimeq}
    \lim_{s\to\infty} \frac{f(e^{s\xi+\sqrt{-1}\nu})}{e^{s(M_1-N_1)(\xi)+ \sqrt{-1} (M_1-N_1)(\nu)}}=1.
  \end{equation}
  Thus,
  \[
    \lim_{s\to\infty} \frac{f(e^{s\xi+\sqrt{-1}\nu})}{|f(e^{s\xi+\sqrt{-1}\nu})|}=\lim_{s\to\infty} \frac{f(e^{s\xi+\sqrt{-1}\nu})}{e^{s(M_1-N_1)(\xi)}}=e^{\sqrt{-1}(M_1-N_1)(\nu)}=\exp(f^t|_C)(e^{\sqrt{-1}\nu}),
  \]
  which gives the second component of (\ref{scalingandtropequation}).

  Now, from (\ref{fraclimeq}) we find
  \[
    \lim_{s\to\infty} \log |f(e^{s\xi+\sqrt{-1}\nu})|-s(M_1-N_1)(\xi)=0.
  \]
Therefore, 
  \begin{align*}
    \lim_{s\to\infty}\frac{1}{s} \log|f(e^{s\xi+\sqrt{-1}\nu})| 
    & =(M_1-N_1)(\xi)\\
    & =f^t(\xi),
  \end{align*}
  which gives the first component of (\ref{scalingandtropequation}).
\end{proof}

\begin{theorem} \label{parttropfunctor}
Partial tropicalization is a functor
\[
PT: \mathbf{PosPoiss}^\bullet_\R\to \mathbf{PTrop}.
\]
\end{theorem}

\begin{proof}
    Recall that after tropicalizing, maps $f:(X,\theta_X,\Phi_X)\to (Y,\theta_Y,\Phi_Y)$ of framed positive varieties with potential preserve the integer lattice and send cones into cones:
    \[
    f^t(\Ccal_{\Phi_X})\subset \Ccal_{\Phi_Y}
    \]
    We then have condition 1 of Definition \ref{partialtrop}, and functoriality will then follow from functoriality of tropicalization, Proposition \ref{tropfunctor}. It remains to check, for a map $f:(X,\pi_X,\theta_X,\Phi_X,\tau_X)\to(Y,\pi_Y,\theta_Y,\Phi_Y,\tau_Y)$, that $PT(f)$ satisfies conditions 2 and 3 of Definition \ref{partialtrop}. 
    
    Extending our notation from above, let $f_s:=(\theta_Y\circ E_{Y,\theta_Y,s})\n\circ f \circ(\theta_X\circ  E_{X,\theta_X,s})$ and let $f_s^\Re:= f_s|_{L_X\times \T_X}$ be the restriction to the real form.
    Then considering the component functions of $f$ and applying Lemma \ref{scalingandtrop}, we have for open linearity chambers $C$ of $f^t$,
    \begin{equation}
    \label{limiteq}
    \lim_{s\to \infty} f_s|_C= f^t|_C\times \exp(f^t|_C):C\times(S^1)^n\to Y^t(\R) \times (S^1)^m,
    \end{equation}
    where $n$ and $m$ are the dimensions of $X$ and $Y$, respectively.
    On each open linearity chamber $C$ of $f^t$, by Proposition \ref{mapsofrealforms} we have in our new coordinates the restricted map
    \begin{equation}
    \label{restmap}
    f_s^\Re |_{C}:(C\cap L_{X})\times \T_{X} \to L_{Y} \times \T_{Y},
    \end{equation}
    has its image in the real form of $Y$. By Lemma \ref{scalingandtrop}, the map $f_s^\Re$ is defined at each point of $C$ for large enough $s$.
    The description of the limit (\ref{limiteq}) gives condition 2. 
    By Proposition \ref{mapsofrealforms}, the map $\theta_Y\n\circ f|_{\Re(X)} \circ \theta_X$ is a Poisson map of real forms when it is defined. Restricting the limit $\lim_{s\to\infty} f_s^{\Re(X)}|_C$ of (\ref{restmap}) to $\Ccal_{\Phi_X}(\R)\times \T_{X}$ and applying
    Theorem \ref{scalingtrans}, we find the map
    \[
    \left. \left( \lim_{s\to\infty} f_s^{\Re(X)}|_{C} \right)\right|_{\Ccal_{\Phi_X}(\R)\times \T_{X}} :(C\cap (\Ccal_{\Phi_X}(\R) \cap L_{X}))\times \T_{X} \to (\Ccal_{\Phi_Y}(\R) \cap L_{Y})\times \T_{Y}.
    \]
    preserves the partially tropicalized Poisson structures. The description of the limit (\ref{limiteq}) gives condition 3.
\end{proof}

\subsection{Partial Tropicalization of \texorpdfstring{$K^*$}{K*}}

Fix a simply-connected semisimple complex Lie group $G$, and let $\mathbf{i}$ be a reduced word for $w_0$. As shown in  Theorem \ref{realformofG^*}, the tuple $(G^*,\sqrt{-1}\pi_{G^*},\bbDelta_\mathbf{i}\n,\Phi_{G^*},\tau)$ is a framed positive Poisson variety with real form, which by Remark \ref{realformG^*remark} is an open dense subset of the Poisson-Lie group $(K^*,\pi_{K^*})$.
We summarize our results in the following theorem.

\begin{theorem} \label{conedescriptiontheorem}
The partial tropicalization $PT(G^*,\sqrt{-1}\pi_{G^*},\bbDelta_\mathbf{i}\n,\Phi_{G^*},\tau)$ is of the form $\Ccal_{\Phi_{BK}}(\R) \times \T$, where $\Ccal_{\Phi_{BK}}$ is the strict extended string cone of Definition \ref{BKpotential}, and $\T$ is a real torus of dimension $d={\rm dim}(N)$. This space is equipped with a constant Poisson bracket and an integrable system.
\end{theorem}

\begin{proof}
    Let $L=L_{G^*,\bbDelta_\mathbf{i}\n,\tau}$. By Theorem \ref{parttropfunctor}, it suffices to describe the cone $L\cap \Ccal_{\Phi_{G^*}}(G^*,\bbDelta_\mathbf{i}\n)(\R)$ for a single reduced word $\mathbf{i}$ of $w_0$. Write $\xi_1,\dots,\xi_{l(w_0)}$ for the coordinates on $(G^*,\bbDelta_\mathbf{i})^t(\R)$ coming from $(\hat{\Delta}_\mathbf{i})_1$; see (\ref{Deltahat})-(\ref{bbDelta}) above. Similarly, write $\xi_1^\tau,\dots,\xi_{l(w_0)}^\tau$ for the coordinates coming from $(\hat{\Delta}_\mathbf{i}\circ \tau)_2$, and $\xi_1^0,\dots,\xi_r^0$ for the coordinates coming from $\Delta_0$. In these coordinates, we have the following expression for $\tau^t:(G^*,\bbDelta_\mathbf{i}\n)^t\to (G^*,\bbDelta_\mathbf{i}\n)^t$:
    \[ 
    \tau^t (\xi_1,\dots,\xi_{l(w_0)},\xi_1^\tau,\dots,\xi_{l(w_0)}^\tau,\xi_1^0,\dots,\xi_r^0)= (\xi_1^\tau,\dots,\xi_{l(w_0)}^\tau,\xi_1,\dots,\xi_{l(w_0)},\xi_1^0,\dots,\xi_r^0).
    \]
    Projecting parallel to the coordinates $\xi_1^\tau,\dots,\xi_{l(w_0)}^\tau$ gives an isomorphism
    \[ 
    L\cong (G^{e,w_0},\Delta_\mathbf{i}\n)^t
    \]
    which takes $\Ccal_{\Phi_{G^*}}$ to $\Ccal_{\Phi_{BK}}$.
    
    A similar analysis in coordinates of $\T=\T_{G^*,\bbDelta_\mathbf{i}\n,\tau} \subset (S^1)^{r+2l(w_0)}$ gives that the dimension of $\T$ is $l(w_0)=\dim(N)$.
\end{proof}

\begin{remark}
Integrable systems of Theorem \ref{conedescriptiontheorem} can be described in more detail. In particular, the constant Poisson bracket $\pi_{PT}$ is of  maximal rank $2d$, symplectic leaves of $\pi_{PT}$ are in one-to-one correspondence with generic symplectic leaves in $K^*$, and symplectic volumes of the leaves of $\pi_{PT}$ coincide with symplectic volumes of the corresponding leaves in $K^*$. This picture will be explored in detail in our forthcoming work. \end{remark}

In particular, we have the following.

\begin{theorem}
For $G={\rm SL}_n(\C)$, the partial tropicalization $(PT(G^*,\sqrt{-1}\pi_{G^*},\bbDelta_\mathbf{i}\n,\Phi_{G^*},\tau),\pi_{PT})$ is isomorphic to the Gelfand-Zeitlin integrable system.
\end{theorem}

\begin{proof}
    In \cite{AD}, this was shown for a specific choice of $\mathbf{i}$. Fix another reduced word $\mathbf{i'}$ for $w_0$. Let 
    \[
    \bbDelta_{\mathbf{i'}}\circ \bbDelta_\mathbf{i}\n: (\C^\times)^n\to (\C^\times)^n
    \]
    be the positive equivalence from the toric chart $\bbDelta_\mathbf{i}\n$ to $\bbDelta_\mathbf{i'}\n$, and let $\{C_j\}_j$ be the set of open linearity chambers of
    \[(\bbDelta_{\mathbf{i'}}\circ \bbDelta_\mathbf{i}\n)^t: (G^*, \bbDelta_\mathbf{i}\n)^t(\R) \to (G^*,\bbDelta_\mathbf{i'}\n)^t(\R).
    \]
    Let
    \[D=\left(G^*, \bbDelta_\mathbf{i}\n)^t(\R)\right)\backslash\cup_j C_j
    \] be the complement of the open linearity chambers $C_j$, and let $\Ccal=\left(\Ccal_{\Phi_{G^*}}(G^*, \bbDelta_\mathbf{i}\n)(\R)\right)\cap L$, where as before we write $L=L_{G^*,\bbDelta_\mathbf{i}\n,\tau}$. By Theorem \ref{parttropfunctor}, the partial tropicalization $PT(\bbDelta_\mathbf{i}\circ \bbDelta_\mathbf{i'}\n)$ determines a Poisson map
    \[
    (PT(G^*,\sqrt{-1}\pi_{G^*},\bbDelta_\mathbf{i}\n,\Phi_{G^*},\tau),\pi_{PT}) \to (PT(G^*,\sqrt{-1}\pi_{G^*},\bbDelta_\mathbf{i'}\n,\Phi_{G^*},\tau),\pi_{PT})
        \]
    away from $D\cap \Ccal \times \T_{G^*,\bbDelta_\mathbf{i}\n,\tau}$, so we just need to show that $\Ccal\not\subset D$. We will show $D\cap \Ccal$ is of positive codimension in $\Ccal$.
        
    First, we note that $D$ is a union of finitely many subspaces $D_k$ of $(G^*,\bbDelta_\mathbf{i}\n)^t(\R)$. So it suffices to show each $D_k\cap \Ccal$ has positive codimension in $\Ccal$.
    
    Second, the cone $\Ccal_{\Phi_{G^*}}(G^*, \bbDelta_\mathbf{i}\n)(\R)$ is open in $(G^*,\bbDelta_\mathbf{i}\n)^t(\R)$, so $\Ccal$ is open in $L$.
    So if some $D_k\cap \Ccal$ has codimension $0$ in $\Ccal$, then $L\subset D_k$.
    
    We then show that $L$ is not contained in any $D_k$. For the reduced word $\mathbf{i}$, consider the coordinates $\xi_1,\dots,\xi_{l(w_0)}$ and $\xi^\tau_1,\dots,\xi^\tau_{l(w_0)}$ introduced in the proof of Theorem \ref{conedescriptiontheorem}. From the definition of $\bbDelta_\mathbf{i}$, we see that $D_k$ is given by some linear equalities of the form
    \[
    \sum_{p=1}^{l(w_0)} a_p \xi_p  =0, \qquad
    \sum_{p=1}^{l(w_0)} b_p \xi_p^\tau =0, \qquad 
    a_p,b_p\in \Q.\]
    On the other hand, $L$ is given by the equations
    \[
    \xi_p=\xi_p^\tau,\qquad p=1,\dots, l(w_0).
    \]
    From this description it is evident that $D_k$ does not contain $L$.
\end{proof}



\begin{thebibliography}{10}

\bibitem{Anton}
A. Alekseev, \emph{On Poisson actions of compact Lie groups on symplectic manifolds}, J. Differential Geom. 45 (1997), no. 2, 241-256.

\bibitem{AD}
A. Alekseev, I. Davydenkova, \emph{Inequalities from Poisson brackets}, Indag. Math. 25 (5) (2014)
846–871.

\bibitem{AM} 
A. Alekseev, E. Meinrenken, \emph{Ginzburg-Weinstein from Gelfand-Zeitlin}, J. Differential Geom.
76 (2007), no. 1, 1-34.

\bibitem{BK1}
A. Berenstein and D. Kazhdan, \emph{Geometric and unipotent crystals}, Geom. Funct. Anal., Special Volume, Part I (2000), 188–236.

\bibitem{BK}
A. Berenstein, D. Kazhdan, \emph{Geometric and unipotent crystals II: From unipotent bicrystals to crystal bases quantum groups}, Contemp. Math., vol. 433, Amer. Math. Soc., Providence, RI, 2007, pp. 13–88.

\bibitem{BFZ}
A. Berenstein, S. Fomin, A. Zelevinsky, \emph{Cluster algebras III. Upper bounds and double Bruhat cells}, Duke Math. J. 126 (1) (2005) 1–52.

\bibitem{BZ1}
A. Berenstein, A. Zelevinsky, \emph{Canonical bases for the quantum group of type $A_r$ and piecewise linear combinatorics}, Duke Math. J. 82 (1996), 473-502.

\bibitem{BZ}
A. Berenstein, A. Zelevinsky, \emph{Total positivity in Schubert varieties}, Comment. Math. Helv. 72 (1997), 128–166. CMP 97:14.

\bibitem{BZ2}
A. Berenstein, A. Zelevinsky: \emph{Tensor product multiplicities, canonical bases and totally positive varieties}. Invent. Math. 143, 77–128 (2001).

\bibitem{ES}
P. Etingof, O. Schiffmann, \emph{Lectures on Quantum Groups}, 2nd edition, International Press, 2002.

\bibitem{FZ}
S. Fomin, A. Zelevinsky, \emph{Double Bruhat Cells and total positivity}, J. Amer. Math. Soc. 12 (1999) 335–380.

\bibitem{GSV}
M. Gekhtman, M. Shapiro, and A. Vainshtein, \emph{Cluster algebras and Poisson geometry}, Mosc. Math. J. 3 (2003), no. 3, 899–934, 1199.

\bibitem{GW}
V. Ginzburg, A. Weinstein, \emph{Lie-Poisson structure on some Poisson Lie groups}, J. Amer. Math. Soc. 5 (2) (1992) 445–453.

\bibitem{KZ}
M. Kogan, A. Zelevinsky, \emph{On symplectic leaves and integrable systems in standard complex semisimple Poisson-Lie groups}, Int. Math. Res. Not. 32 (2002) 1685–1702.

\bibitem{LW}
J. H. Lu, A. Weinstein, \emph{Poisson-Lie groups, dressing transformations and Bruhat decompositions}, J. Differential Geom. 31 (1990), no.2, 501–526.

\bibitem{Mi}
G. Mikhalkin, \emph{Tropical geometry and its applications}, International Congress of Mathematicians. Vol. II, 827–852, Eur. Math. Soc., Zürich, 2006.

\bibitem{Ri}
K. Rietsch, \emph{A mirror symmetric construction of $qH^{*}_T(G/P)_{(q)}$}, Adv. Math. 217 (2008), no. 6, 2401–2442. 

\bibitem{STS}
M. A. Semenov-Tian-Shansky, \emph{What is a classical $r$-matrix?}, Functional
Analysis and Its Applications 17 (1983) no. 4, 259–272.

\bibitem{STS2}
M. A. Semenov-Tian-Shansky, \emph{Dressing transformations and Poisson group actions}, Publ. Res. Inst. Math. Sci. 21 (1985), no. 6, 1237–1260.

\bibitem{Xu}
P. Xu, \emph{Dirac submanifolds and Poisson involutions}, Ann. Scient. \'{E}c. Norm. Sup., 36 (2003), 403 - 430

\end{thebibliography}
\end{document}